\documentclass[twocolumn]{autart} 

\usepackage{amsmath}
 \usepackage{subfig}

\usepackage{graphicx}
\usepackage{amsfonts}
\usepackage{epstopdf}
\newtheorem{theorem}{Theorem}
\usepackage{color}
\usepackage{xcolor}
\usepackage[toc,page]{appendix}
\usepackage[english]{babel} 
\usepackage{amssymb}
\usepackage{amsmath}
\usepackage{mathdots}

\begin{document}

\setlength\abovedisplayskip{7pt}
\setlength\belowdisplayskip{7pt}
\setlength\abovedisplayshortskip{7pt}
\setlength\belowdisplayshortskip{7pt}

\begin{frontmatter}

\title{Event-Triggered Control of a Continuum Model of Highly Re-Entrant Manufacturing System} 

\thanks[footnoteinfo]{Corresponding author: M.~Diagne. Tel. +001-518-276-8145. }

\author[UMICH]{Mamadou Diagne}\ead{diagnm@rpi.edu} \, and \, 
\author[Karafyllis]{Iasson Karafyllis}\ead{iasonkar@central.ntua.gr}

\address[UMICH]{Department of Mechanical Aerospace and Nuclear Engineering, Rensselaer Polytechnic Institute, Troy, New York, 12180, USA} 
\address[Karafyllis]{Department of Mathematics, National Technical University of Athens,
Zografou Campus, 15780, Athens, Greece}

\begin{abstract}                          
 With the unceasing growth of intelligent production lines that integrate sensors, actuators, and controllers in a wireless communication environment via internet of things (IoT), we design an event-triggered boundary controller for a continuum model of highly re-entrant manufacturing systems for which the \textcolor{black}{influx rate} of product is the controlled quantity. The designed controller can potentially operate in networked control systems subject to limited information sharing resources. A Lyapunov argument is utilized to derive the boundary controller together with a feasible event generator that avoids the occurrence of Zeno behavior for the closed-loop system. The global stability estimate is established using the logarithmic norm of the state due to the system's nonlinearity \textcolor{black}{and positivity of the density}. Furthermore,  robustness of the proposed controller \textcolor{black}{with respect to} the sampling schedule and sampled-data \textcolor{black}{stabilization results} are established. Consistent simulation results that support the proposed theoretical statements are provided.
\end{abstract}

\end{frontmatter}


\section{Introduction}\label{s1}

The comprehensive development of event-triggered control in the past few years has been motivated by the crucial need to preserve limited computational and communication resources during the execution of feedback control tasks in constrained Networked Control Systems (NCS). The key idea of event-triggered control consists of the starting of ``only necessary" control action when events generated by the real-time systems' response occur, reducing considerably the number of execution of control tasks while preserving satisfactory closed-loop system performance. In comparison to the well-known periodical sampled-data approach, event-triggered enables noncyclic updates of control signals and thereby offers more flexibility to manage constraints arising from networks and systems interactions. Generally speaking, event-based control designs require two challenging steps that  can be achieved simultaneously or independently: 
\begin{itemize}
\item The construction of a feasible successive execution time sequence (with positive minimum dwell-time), \textcolor{black}{which is related} to the event generator that determines the time instants at which the control input is updated. 
\item The design of the feedback control signal that ensures closed-loop performance specifications.
\end{itemize}
The event-based PID control design proposed in \cite{s9} and the event-based sampling for first-order stochastic systems \cite{s8} serve as pioneering contributions in the field. These results consider scheduling algorithms that only recondition feedback control signals according to an error with respect to a given state norm. Later on, motivated by the stabilization of constrained  networked control systems \cite{s4}, major results dedicated to both linear \cite{s1,s3,s5,s6} and  nonlinear \cite{s2,s7,s10,s11,s12} finite-dimensional systems have been established introducing Lyapunov-based triggering conditions to ensure stability at desired decay rate. As well, advanced control designs involving robustness analysis \cite{s13,s14}, stabilization of multi-agent consensus problems and adaptive control \cite{s16} can be found in the abundant literature. From a practical point of view, real implementation of event-based schemes has been achieved to control wireless throttling valves \cite{s34}, the angular position of a gyroscope \cite{s35}, the formation of a group of VTOL-UAVs \cite{s37} and a greenhouse temperature \cite{s36}, to name a few.

Concurrently, substantial efforts have been taken to develop event-based control of infinite-dimensional systems. In this case, the deduction of the triggering condition follows a Lyapunov criterion and is analogous to that of finite-dimensional systems. 
Early results are developed upon reduced-order models that describe the dominant dynamics of reaction-diffusion systems. In this case, \textcolor{black}{the resulting} linear finite-dimensional systems are exploited to match the control objectives \cite{s31,s21}. However, it is well known that the order of approximation is not trivially determined a priori using a modal expansion of PDEs (Partial Differential Equations). Alongside, both event-triggered and sampled-data control have been successfully developed for reaction-advection-diffusion PDE \cite{s17,s18,s20,s27,s32,s25a} and ODE-PDE cascading systems \cite{s30} without model reduction. \textcolor{black}{We emphasize that \cite{s25a}, which employs a small gain design, is one of the first attempt of  event-triggered boundary control of 1D parabolic PDEs.} For hyperbolic PDEs substantial developments can be found in \cite{s19,s22,s23,s24,s25,s26,s28}. It is worth to mention that local stability result has been recently achieved applying sampled-data control to a nonlinear PDE governed by 1-D Kuramoto-Sivashinsky equation \cite{s29}.

  Furthermore, as emphasized in the premier contribution \cite{s9}, for factory lines, the event-based nature of the sampling can be related to the process's intrinsic production rate. From a modeling point of view, discrete-event \cite{s47,s48}, Effective Processing Time (EPT) \cite{s51} and clearing functions \cite{s50} representations, which requires a deep knowledge of various process specifications or only accounts on arrival and departure events of the parts to a workstation when Work in Progress (WIP) is not strongly varying, has been proven inefficient for operation planning and control. Motivated by the study of the transient behavior of manufacturing systems with high fluctuations of the WIP, nonlocal transport PDE models \cite{s43,s44,s49} has emerged during the past few years. These continuum models describe the time evolution of the flow of manufactured products using the spatial distribution of product density as a key variable. Several contributions considering the control of boundary influx of parts with PI controller \cite{s45}, Lyapunov-based design \cite{s40,s41,s42}, \textcolor{black}{small gain design \cite{s33}}, predictor-feedback design \cite{s46} or optimal \cite{s38,s39} control techniques of these non-local PDEs have been recently developed by researchers.


With the unceasing growth of intelligent production lines that integrate sensors, actuators, and controllers in a wireless communication environment via internet of things (IoT), we design an event-triggered boundary controller for a continuum model of highly re-entrant manufacturing systems for which the \textcolor{black}{influx rate} of product is the controlled quantity. The designed controller can potentially operate in networked control systems subject to limited information sharing resources. A Lyapunov argument is utilized to derive the boundary controller together with a feasible event generator that avoids the occurrence of Zeno behavior for the closed-loop system. The global stability estimate is established using the logarithmic norm of the state due to the system's nonlinearity \textcolor{black}{and positivity of the density}. Furthermore,  robustness of the proposed controller \textcolor{black}{with respect to} the sampling schedule and sampled- data \textcolor{black}{stabilization results (Theorem \ref{T4})}are established. Consistent simulation results that support the proposed theoretical statements are provided. 

The paper is organized as follows: Section \ref{s2} describes the modeling of highly re-entrant manufacturing systems as nonlocal PDE with an influx boundary condition. The existence and uniqueness of solutions for a piecewise continuous control input signal are stated in Section \ref{s3}. The construction of robust event-triggered and sampled-data stabilizing boundary controllers is discussed in Section \ref{s4}. Simulation results that demonstrate the feasibility of both event-triggered and sampled-data controllers are discussed in Section \ref{s5}. Finally, the paper ends with concluding remarks and future research directions in Section \ref{s6}.

     \bigskip 
 \noindent

\noindent \textbf{Notation:} Throughout the paper, we adopt the following notation.

\begin{itemize} \item  $\mathbb{R} _{+} :=[0,+\infty )$. Let $u:\mathbb{R} _{+} \times [0,1]\to \mathbb{R} $ be given. We use the notation $u[t]$ to denote the profile at certain $t\ge 0$, i.e., $(u[t])(x)=u(t,x)$ for all $x\in [0,1]$. For bounded functions $f:I\to \mathbb{R} $, where $I\subseteq \mathbb{R} $ is an interval, we set $\left\| f\right\| _{\infty } ={\mathop{\sup }\limits_{x\in I}} \left(\left|f(x)\right|\right)<+\infty $. We use the notation $f'(x)$ for the derivative at $x\in [0,1]$ of a differentiable function $f:[0,1]\to \mathbb{R} $. 

 \item Let $S\subseteq \mathbb{R} ^{n} $ be an open set and let $A\subseteq \mathbb{R} ^{n} $ be a set that satisfies $S\subseteq A\subseteq cl(S)$. By $C^{0} (A\; ;\; \Omega )$, we denote the class of continuous functions on $A$, which take values in $\Omega \subseteq \mathbb{R} ^{m} $. By $C^{k} (A\; ;\; \Omega )$, where $k\ge 1$ is an integer, we denote the class of functions on $A\subseteq \mathbb{R} ^{n} $, which takes values in $\Omega \subseteq \mathbb{R} ^{m} $ and has continuous derivatives of order $k$. In other words, the functions of class $C^{k} (A;\Omega )$ are the functions which have continuous derivatives of order $k$ in $S=int(A)$ that can be continued continuously to all points in $\partial S\cap A$.  When $\Omega =\mathbb{R} $ then we write $C^{0} (A\; )$ or $C^{k} (A\; )$. 
 \item A left-continuous function $f:[0,1]\to \mathbb{R} $ (i.e. a function with ${\mathop{\lim }\limits_{y\to x^{-} }} \left(f(y)\right)=f(x)$ for all $x\in (0,1]$) is called piecewise $C^{1} $ on $[0,1]$ and we write $f\in PC^{1} ([0,1])$, if the following properties hold: (i) for every $x\in [0,1)$ the limits ${\mathop{\lim }\limits_{y\to x^{+} }} \left(f(y)\right)$, ${\mathop{\lim }\limits_{h\to 0^{+} ,y\to x^{+} }} \left(h^{-1} \left(f(y+h)-f(y)\right)\right)$ exist and are finite, (ii) for every $x\in (0,1]$ the limit ${\mathop{\lim }\limits_{h\to 0^{-} }} \left(h^{-1} \left(f(x+h)-f(x)\right)\right)$ exists and is finite, (iii) there exists a set $J\subset (0,1)$ of finite cardinality, where $f'(x)={\mathop{\lim }\limits_{h\to 0^{-} }} \left(h^{-1} \left(f(x+h)-f(x)\right)\right)={\mathop{\lim }\limits_{h\to 0^{+} }} \left(h^{-1} \left(f(x+h)-f(x)\right)\right)$ holds for $x\in (0,1)\backslash J$, and (iv) the mapping $((0,1)\backslash J)\ni x\to f'(x)\in \mathbb{R} $ is continuous. Notice that we require a piecewise $C^{1} $ function to be left-continuous but not continuous. 
\end{itemize}
\section{Continuum Model of Highly re-entrant Manufacturing Systems}\label{s2}
 Manufacturing systems with a high volume and a large number of consecutive production steps (which typically number in the many hundreds) are often modeled by non-local PDEs \cite{s43,s44,s49}. The general one-dimensional continuity equation expressing mass conservation along the  production stages is considered. Defining the flow of unit parts  per unit time as $\textcolor{black}{(F(\rho[t]))(x) }$ where $F$ is  the flux at the production stage $x$ and time $t$ depending of  the density of product $\rho(t,x)$, namely, the work in progress, the following  equation of conservation can be written between two production stages $x_1$ and $x_2$
 \begin{equation} \label{C0} 
\frac{\partial}{\partial \, t} \int_{x_1}^{x_2}\rho(t,x)dx= \textcolor{black}{(F(\rho[t]))(x_1) }-\textcolor{black}{(F(\rho[t]))(x_2) }, 
\end{equation} 
and equivalently in  \textcolor{black}{differential} form as 
\begin{equation} \label{C1} 
\frac{\partial \, \rho (t,x)}{\partial \, t} +\frac{\partial  \textcolor{black}{(F(\rho[t]))(x) } }{\partial \, x} =0. 
\end{equation} 
The density function  $\rho (t,x)$  defined at time $t\ge 0$ and  stage $x\in [0,1]$, is  required to be positive (i.e., $\rho (t,x)>0$ for $(t,x)\in \mathbb{R} _{+} \times [0,1]$). Moreover, it  has  a spatially uniform and positive equilibrium profile $$\rho (x)\equiv \rho _{s}, \quad  \rho _{s} >0. $$
 Controlling of the governing equation \eqref{C1} consists of determining the influx 
\begin{equation} \label{C2} 
\textcolor{black}{(F(\rho[t]))(0)}=u(t) 
\end{equation} 
that results in  the desired outflux
\begin{equation} \label{C3} 
(F(\rho[t]))(1)=y(t).
\end{equation} 
A challenging aspect of the model \eqref{C1} is the characterization of  flux  $F$ as a function of   the work in progress $\rho$. Following \cite{s43,s44,s49}, $F$ can be  defined as a function of the length of the queue or the total work load  $W(t)$, that is, 
$$\textcolor{black}{(F(\rho[t]))(x)}=\lambda \left(W(t)\right)\rho(t,x),$$ leading to the following equation
\begin{align} \label{GrindEQ__1_} 
\frac{\partial \, \rho }{\partial \, t} (t,x)+\lambda \left(W(t)\right)\frac{\partial \, \rho }{\partial \, x} (t,x)=0, 
\end{align} 
where 
\begin{align} \label{GrindEQ__2_} 
W(t)=\int _{0}^{1}\rho (t,x)dx. 
\end{align} 
Here,  $\lambda \in C^{1} \left(\mathbb{R} _{+} ;(0,+\infty )\right)$ is a non-increasing function that determines the production speed.  Hence, the process  influx rate at the boundary defined in  \eqref{C2} becomes
\begin{align} \label{GrindEQ__3_} 
\rho (t,0)\lambda \left(W(t)\right)=u(t) 
\end{align} 
where $u(t)\in (0,+\infty )$ is the control input.

\begin{figure}[htbp]
\centering
\includegraphics[width=0.45\textwidth]{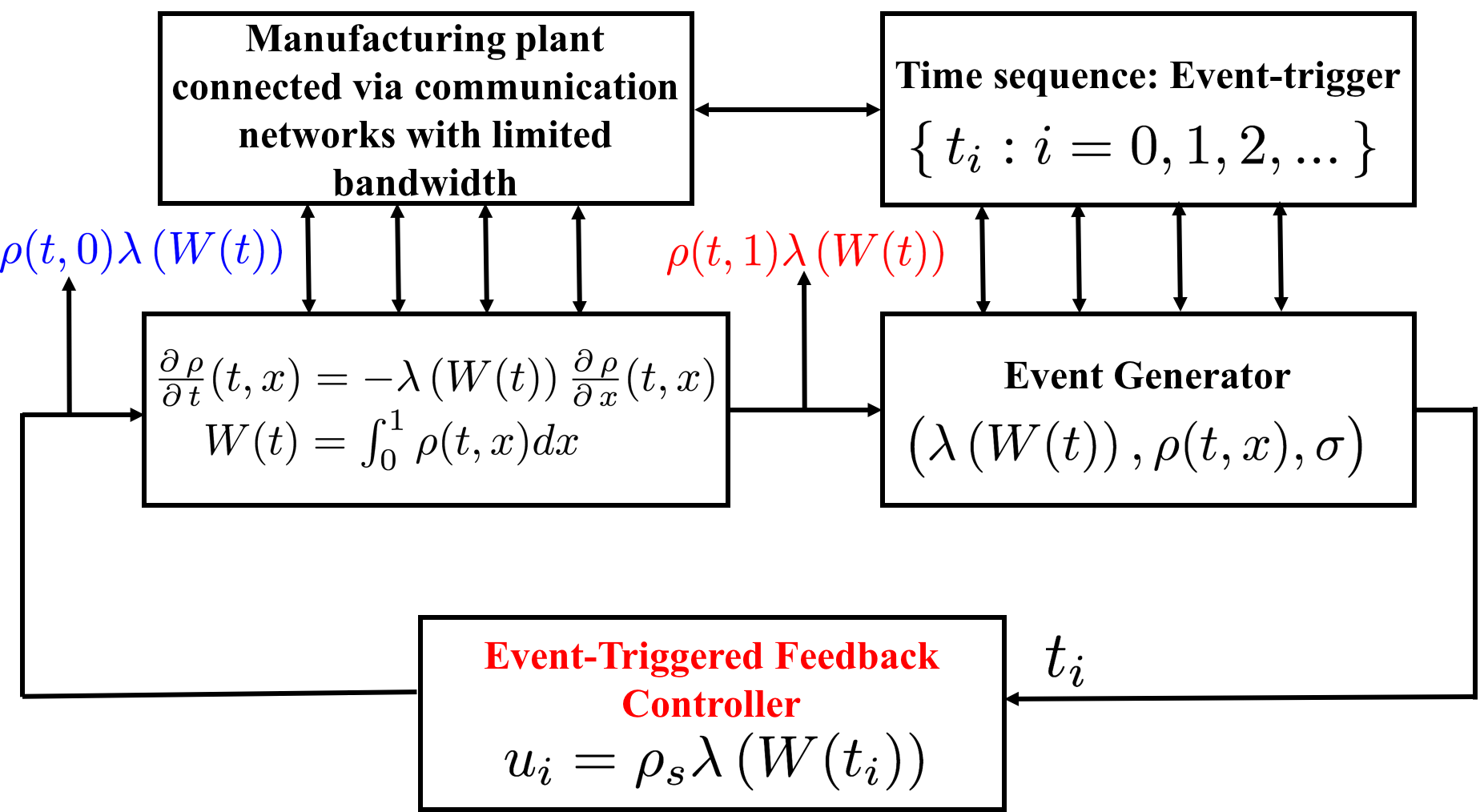}
\caption{Event-trigger closed-loop system.}\label{system_figure-struct}
\end{figure}

In what follows, we intend to apply event-triggered boundary control through $u(t)\in (0,+\infty )$ to the  the manufacturing plant described by the \eqref{GrindEQ__1_}--\eqref{GrindEQ__3_}, which will achieve global stabilization of the spatially uniform equilibrium profile $\rho (x)\equiv \rho _{s} $. Therefore, the control $u(t)\in (0,+\infty )$ will have the form
\begin{align} \label{GrindEQ__4_}
u(t)=u_{i} , \ for \ t\in [t_{i} ,t_{i+1} )                                                          
\end{align}
where $u_{i} >0$ are the input values that will be determined by the controller and the times $\left\{\, t_{i} \, :\, i=0,1,2,...\, \right\}$ will be the times of the events that will be determined by the event-trigger, which will constitute an increasing sequence with $t_{0} =0$. The structure of the closed-loop system consisting of the plant, the event generator and the control input is depicted in Figure \ref{system_figure-struct}.

\section{ Notion of Solution} \label{s3}

First we  describe precisely  the notion of the solution of the closed-loop system. The following result plays an instrumental role in the construction of solutions for the closed-loop system \eqref{GrindEQ__1_}, \eqref{GrindEQ__2_}, \eqref{GrindEQ__3_}, \eqref{GrindEQ__4_}. \\

\begin{theorem}\label{T1} Consider the initial-boundary value problem \eqref{GrindEQ__1_}, \eqref{GrindEQ__2_}, \eqref{GrindEQ__3_} with initial condition
\begin{align} \label{GrindEQ__5_}
\rho (0,x)=\rho _{0} (x), \ for \ x\in (0,1]   
\end{align}
where $\lambda \in C^{1} \left(\mathbb{R} _{+} ;(0,+\infty )\right)$ is a non-increasing function, $\rho _{0} \in PC^{1} \left([0,1]\right)$, $u(t)\equiv u>0$ and ${\mathop{\inf }\limits_{x\in (0,1]}} \left(\rho _{0} (x)\right)>0$. Suppose that there exists a constant $K>0$ such that $\left|\lambda '(s)\right|\le K$ for all $s\ge 0$. Then there exists $t_{\max } \in (0,+\infty ]$ and unique functions $W\in C^{0} \left([0,t_{\max } );(0,+\infty )\right)$, $\rho :[0,t_{\max } )\times [0,1]\to (0,+\infty )$ with $\rho [t]\in PC^{1} \left([0,1]\right)$ and ${\mathop{\inf }\limits_{x\in [0,1]}} \left(\rho (t,x)\right)>0$ for all $t\in [0,t_{\max } )$, which constitute a solution of the initial-boundary value problem \eqref{GrindEQ__1_}, \eqref{GrindEQ__2_}, \eqref{GrindEQ__3_}, \eqref{GrindEQ__5_} in the following sense:
\begin{itemize}
\item The function  $\rho:[0,t_{\max } )\times [0,1]\to (0,+\infty )$ is of class $C^{1} \left([0,t_{\max } )\times [0,1]\backslash \Omega \right)$, where 
\begin{align} \label{GrindEQ__6_a} 
\Omega =\Omega_1\cup\Omega_2 \end{align} 
with
\begin{align} \label{GrindEQ__6_b} 
\hspace{-0.5cm}\Omega_1 &={\mathop{\cup }\limits_{i=0,...,N}} \left\{ (t,r_{i} (t)):t\in [0,t_{\max } ),r_{i} (t)\le 1 \right\},\\
 \label{GrindEQ__6_c} 
\Omega_2 &={\mathop{\cup }\limits_{i=0,...,N}} \Bigg\{\, (t,r_{i} (t)-1):\, \, t\in [0,t_{\max } ),\nonumber \\ &~~~~~~~~~~~~~~~~~~~1\le r_{i} (t)\le 2\, \Bigg\}
\end{align} 
 and \begin{align}r_{i} (t)=\xi _{i} +\int _{0}^{t}\lambda (W(s))ds, \quad  (i=0,...,N),\end{align} $\xi _{0} =0$ and $\xi _{i} \in [0,1)$ ($i=1,...,N$) are the points (in increasing order) for which $\rho _{0} \in C^{1} \left([0,1]\backslash \{ \xi _{0} ,...,\xi _{N} \} \right)$, 

\item Equation \eqref{GrindEQ__1_} holds for all $(t,x)\in [0,t_{\max } )\times [0,1]\backslash \Omega $ and equations \eqref{GrindEQ__2_}, \eqref{GrindEQ__3_}, \eqref{GrindEQ__5_} hold for $t\in [0,t_{\max } )$ with $u(t)\equiv u>0$.

\item  Finally, if $t_{\max } <+\infty $ then \begin{align}{\mathop{\lim \sup }\limits_{t\to t_{\max }^{-} }} \left(\left\| \rho [t]\right\| _{\infty } \right)=+\infty. \end{align} \end{itemize} \end{theorem}
\textcolor{black}{ The detailed  proof  of Theorem \ref{T1} is given in Section Appendix}.

 Given $\rho _{0} \in PC^{1} \left([0,1]\right)$ with ${\mathop{\inf }\limits_{x\in (0,1]}} \left(\rho _{0} (x)\right)>0$ and using Theorem \ref{T1}, we are in a position to construct a unique solution for the closed-loop system \eqref{GrindEQ__1_}, \eqref{GrindEQ__2_}, \eqref{GrindEQ__3_}, \eqref{GrindEQ__4_}, \eqref{GrindEQ__5_} by means of the following algorithm for each integer $i\ge 0$:
 
\begin{enumerate}
\item  Given $t_{i} \ge 0$ and $\rho [t_{i} ]\in PC^{1} \left([0,1]\right)$ with ${\mathop{\inf }\limits_{x\in (0,1]}} \left(\rho (t_{i} ,x)\right)>0$, determine $u_{i} >0$ by using an appropriate feedback law.
\item  Construct the solution of the initial-boundary value problem 
\begin{align*}&\frac{\partial \, y}{\partial \, t} (t,x)+\lambda \left(V(t)\right)\frac{\partial \, y}{\partial \, x} (t,x)=0\\
&V(t)=\int _{0}^{1}y(t,x)dx \\ 
&y(t,0)\lambda \left(V(t)\right)=u_{i}, \end{align*} with initial condition $y(0,x)=\rho (t_{i} ,x)$, for $x\in (0,1]$. Theorem 1 guarantees that such a solution exists as long as $y$ remains bounded in the interval $[0,t_{i+1} -t_{i} ]$.
\item  Using the event-trigger and assuming that $\rho [t_{i} +s]=y[s]$ for $s\ge 0$, determine the time of the next event $t_{i+1} >t_{i} $. 

\item  Set $\rho [t_{i} +s]=y[s]$ for $0\le s\le t_{i+1} -t_{i} $ and repeat. 
\end{enumerate}

\noindent The above algorithm guarantees that the solution satisfies $\rho [t]\in PC^{1} \left([0,1]\right)$ and ${\mathop{\inf }\limits_{x\in [0,1]}} \left(\rho (t,x)\right)>0$ for all $t\in [0,t_{\max } )$, where $t_{\max } \in (0,+\infty ]$ is the maximal existence time of the solution. Moreover, equations \eqref{GrindEQ__2_}, \eqref{GrindEQ__3_}, \eqref{GrindEQ__4_}, \eqref{GrindEQ__5_} hold, while the PDE \eqref{GrindEQ__1_} holds for $(t,x)\in [0,t_{\max } )\times [0,1]$ almost everywhere. Finally, if $t_{\max } <{\mathop{\lim }\limits_{i\to +\infty }} \left(t_{i} \right)$ then ${\mathop{\lim \sup }\limits_{t\to t_{\max }^{-} }} \left(\left\| \rho [t]\right\| _{\infty } \right)=+\infty $, which implies that the solution exists as long as it is bounded from above.
\section{ Event-triggered  and Sampled-data Control design}\label{s4}
\subsection{Event-Triggered Control}
 Let $\sigma >0$ be a parameter of the controller and consider the event-triggered controller given by the following formulas for all integers $i\ge 0$:
\begin{align} \label{GrindEQ__7_} 
t_{0} =0 
\end{align} 
\begin{align} \label{GrindEQ__8_a} 
\mu_{i} =\inf \bigg\{\tau & >t_{i}: \left|\ln \left(\frac{\rho (\tau ,0)}{\rho _{s} } \right)\right|>\nonumber \\ &\exp \left(-\sigma \int _{t_{i} }^{\tau }\lambda (W(s))ds \right)\nonumber\\ &{\mathop{\sup }\limits_{0<x\le 1}} \left(\left|\ln \left(\frac{\rho (t_{i} ,x)}{\rho _{s} } \right)\right|\exp (-\sigma x)\right)\, \bigg\} \\ \label{GrindEQ__8_b}t_{i+1} =\min& \left ( t_{i} +\frac{1}{\lambda (0)} \, ,\, \mu_{i} \, \right)
\end{align} 
\begin{align} \label{GrindEQ__9_} 
u_{i} =\rho _{s} \lambda \left(W(t_{i} )\right) 
\end{align} 
 It is possible to show that the event-triggered controller \eqref{GrindEQ__4_}, \eqref{GrindEQ__7_}, \eqref{GrindEQ__8_a}, \eqref{GrindEQ__8_b}, \eqref{GrindEQ__9_} guarantees global stabilization of system \eqref{GrindEQ__1_}, \eqref{GrindEQ__2_}, \eqref{GrindEQ__3_}. Moreover, no Zeno behavior can occur for the closed-loop system. Therefore, we state the following result. 

\begin{theorem}\label{T2}Suppose that there exists a constant $K>0$ such that $\left|\lambda '(s)\right|\le K$ for all $s\ge 0$. Let $\sigma >0$ be a given parameter and consider the closed-loop system \eqref{GrindEQ__1_}, \eqref{GrindEQ__2_}, \eqref{GrindEQ__3_}, \eqref{GrindEQ__4_}, \eqref{GrindEQ__7_}, \eqref{GrindEQ__8_a}, \eqref{GrindEQ__8_b}, \eqref{GrindEQ__9_}. Then there exists a non-increasing function $T:\mathbb{R} _{+} \to (0,+\infty )$ such that for every $\rho _{0} \in PC^{1} \left([0,1]\right)$ with ${\mathop{\inf }\limits_{x\in (0,1]}} \left(\rho _{0} (x)\right)>0$, the following estimates hold for the solution $\rho :\mathbb{R} _{+} \times [0,1]\to (0,+\infty )$ with $\rho [t]\in PC^{1} \left([0,1]\right)$ and ${\mathop{\inf }\limits_{x\in [0,1]}} \left(\rho (t,x)\right)>0$ for all $t\ge 0$ of the initial-boundary value problem \eqref{GrindEQ__1_}, \eqref{GrindEQ__2_}, \eqref{GrindEQ__3_}, \eqref{GrindEQ__4_}, \eqref{GrindEQ__5_}, \eqref{GrindEQ__7_}, \eqref{GrindEQ__8_a}, \eqref{GrindEQ__8_b}, \eqref{GrindEQ__9_}:
\begin{align} \label{GrindEQ__10_}
{\mathop{\sup }\limits_{0\le x\le 1}}& \left(\left|\ln \left(\frac{\rho (t,x)}{\rho _{s} } \right)\right|\right)\le \exp \left(-\sigma (c(\rho _{0} )t-1)\right)\nonumber  \\& \times {\mathop{\sup }\limits_{0<x\le 1}} \left(\left|\ln \left(\frac{\rho _{0} (x)}{\rho _{s} } \right)\right|\right), \ \textrm{for} \ \textrm{all} \ \ t \ge 0                    
\end{align} 
\begin{align} \label{GrindEQ__11_}
t_{i+1} \ge t_{i} +T\left({\mathop{\sup }\limits_{0<x\le 1}} \left(\left|\ln \left(\frac{\rho _{0} (x)}{\rho _{s} } \right)\right|\right)\right),\ \textrm{for} \ \textrm{all} \ \ t \ge 0                                        
\end{align} 
where 
\begin{align} \label{GrindEQ__12_} 
c(\rho _{0} ):=\lambda \left(\rho _{s} \exp \left(\exp (\sigma ){\mathop{\sup }\limits_{0<x\le 1}} \left(\left|\ln \left(\frac{\rho _{0} (x)}{\rho _{s} } \right)\right|\right)\right)\right) 
\end{align} 
\end{theorem}
\textbf{Proof of Theorem \ref{T2}:}
 Let  $\rho _{0} \in PC^{1} \left([0,1]\right)$ with ${\mathop{\inf }\limits_{x\in (0,1]}} \left(\rho _{0} (x)\right)>0$ be given (arbitrary). \textcolor{black}{By virtue of Theorem 1}, there exists a unique solution $\rho:\left[0,t_{\max } \right)\times [0,1]\to (0,+\infty )$ with $\rho [t]\in PC^{1} \left([0,1]\right)$ and ${\mathop{\inf }\limits_{x\in [0,1]}} \left(\rho (t,x)\right)>0$, for all $\ t\in \left[0,t_{\max } \right)$, of the initial-boundary value problem \eqref{GrindEQ__1_}, \eqref{GrindEQ__2_}, \eqref{GrindEQ__3_}, \eqref{GrindEQ__4_}, \eqref{GrindEQ__5_}, \eqref{GrindEQ__7_}, \eqref{GrindEQ__8_a}, \eqref{GrindEQ__8_b}, \eqref{GrindEQ__9_}. The solution $\rho :\left[0,t_{\max } \right)\times [0,1]\to (0,+\infty )$ satisfies for all $i\ge 0$ with $t_{i} <t_{\max } $ and $t\in \left[t_{i} ,\min (t_{i+1} ,t_{\max } )\right)$:
\begin{align} \label{GrindEQ__18_} 
\rho (t,x)=\left\{\begin{array}{c} {\rho \left(t_{i} ,x-\int _{t_{i} }^{t}v(s)ds \right)\;  if\,\, 1\ge x>\int _{t_{i} }^{t}v(s)ds } \\ {\frac{u_{i} }{v\left(\tilde{t}(t,x)\right)} \,\,  if\,\, 0\le x\le \int _{t_{i} }^{t}v(s)ds } \end{array}\right.  
\end{align} 
 where $\tilde{t}(t,x)\in [t_{i} ,t]$ is the unique solution of the equation 
 \begin{align}x&=\int _{\tilde{t}(t,x)}^{t}v(s)ds,\\ 
 v(t)&=\lambda \left(W(t)\right).
 \end{align} 
 Formula \eqref{GrindEQ__18_} is a consequence of \eqref{GrindEQ__1_}, \eqref{GrindEQ__2_}, \eqref{GrindEQ__3_} and \eqref{GrindEQ__4_}. Notice that since $v(t)\le \lambda (0)$ (a direct consequence of the fact that $\lambda \in C^{1} \left(\mathbb{R} _{+} ;(0,+\infty )\right)$ is a non-increasing function) and $t_{i+1} -t_{i} \le \frac{1}{\lambda (0)} $ (a direct consequence of (8)), it follows that \textcolor{black}{formula \eqref{GrindEQ__18_} } is valid for all $t\in \left[t_{i} ,\min (t_{i+1} ,t_{\max } )\right)$. 

 We next define for all $t\in \left[0,t_{\max } \right)$:
\begin{align} \label{GrindEQ__19_} 
V(t)={\mathop{\sup }\limits_{0\le x\le 1}} \left(\left|\ln \left(\frac{\rho (t,x)}{\rho _{s} } \right)\right|\exp (-\sigma x)\right) 
\end{align} 
 Combining \eqref{GrindEQ__18_} and \eqref{GrindEQ__19_} we get for all $i\ge 0$ with $t_{i} <t_{\max } $ and $t\in \left[t_{i} ,\min (t_{i+1} ,t_{\max } )\right)$:
\begin{align} \label{GrindEQ__20_} 
 \max &\Bigg\{\sup \limits_{0<\xi \le 1-\int _{t_{i} }^{t}v(s)ds } \Bigg(\left|\ln \left(\frac{\rho \left(t_{i} ,\xi \right)}{\rho _{s} } \right)\right|\nonumber\\ &\times \exp \left(-\sigma \left(\int _{t_{i} }^{t}v(s)ds +\xi \right)\right)\Bigg),\nonumber\\
& \sup \limits_{0\le x\le \int _{t_{i} }^{t}v(s)ds }\Bigg( \left|\ln \left(\frac{\rho \left(\tilde{t}(t,x),0\right)}{\rho _{s} } \right)\right| \exp (-\sigma x)\Bigg)\Bigg\}\nonumber\\
&\leq \max \Bigg\{\exp \left(-\sigma \int _{t_{i} }^{t}v(s)ds \right)V(t_{i} ),\nonumber \\
&\sup \limits_{t_{i} \le \tau \le t} \left(\left|\ln \left(\frac{\rho \left(\tau ,0\right)}{\rho _{s} } \right)\right|\exp \left(-\sigma \int _{\tau }^{t}v(s)ds \right) \right)\Bigg\}
\end{align} 
For the derivation of \eqref{GrindEQ__20_}, we have used the fact that $\tilde{t}(t,x)\in [t_{i} ,t]$ is continuously decreasing with respect to $x$ with $\tilde{t}(t,0)=t$ and $\tilde{t}\left(t,\int _{t_{i} }^{t}v(s)ds \right)=t_{i} $.  

Combined with definition \eqref{GrindEQ__19_}, the event-trigger \eqref{GrindEQ__8_a}, \eqref{GrindEQ__8_b}, (as well as continuity of $v(t)=\lambda \left(W(t)\right)$, which implies continuity of $\rho (t,0)=\frac{u_{i} }{v(t)} $ for $t\in \left[t_{i} ,\min (t_{i+1} ,t_{\max } )\right)$) gives for all $i\ge 0$ with $t_{i} <t_{\max } $:
\begin{align} \label{GrindEQ__21_}
\left|\ln \left(\frac{\rho (\tau ,0)}{\rho _{s} } \right)\right|\le \exp \left(-\sigma \int _{t_{i} }^{\tau }\lambda (W(s))ds \right)V(t_{i} ),                      
\end{align} 
for  all $\tau \in \left[t_{i} ,\min (t_{i+1} ,t_{\max } )\right)$.

 Combining \eqref{GrindEQ__20_} and \eqref{GrindEQ__21_}, we get for all $i\ge 0$ with  $t\in \left[t_{i} ,\min (t_{i+1} ,t_{\max } )\right)$ and $t_{i} <t_{\max } $: 
\begin{align} \label{GrindEQ__22_} 
V(t)\le \exp \left(-\sigma \int _{t_{i} }^{t}v(s)ds \right)V(t_{i} ).
\end{align} 
 Definition \eqref{GrindEQ__19_} and inequality \eqref{GrindEQ__22_} show that $\left\| \rho [t]\right\| _{\infty } $ is bounded on $\left[t_{i} ,\min (t_{i+1} ,t_{\max } )\right)$. Consequently, Theorem \ref{T1} guarantees that $t_{i+1} <t_{\max } $. Therefore, $t_{\max } ={\mathop{\lim }\limits_{i\to +\infty }} \left(t_{i} \right)$. 

Notice that the notion of solution that we have adopted guarantees that \eqref{GrindEQ__22_} is also valid for $t=t_{i+1} $. Applying \eqref{GrindEQ__22_} inductively, we get for all $t\in \left[0,{\mathop{\lim }\limits_{i\to +\infty }} \left(t_{i} \right)\right)$: 
\begin{align} \label{GrindEQ__23_} 
V(t)\le \exp \left(-\sigma \int _{0}^{t}v(s)ds \right)V(0).
\end{align} 
Since $v(t)=\lambda(W(t))\geq 0$ for all $t\in \left[0,{\mathop{\lim }\limits_{i\to +\infty }} \left(t_{i} \right)\right)$, from \eqref{GrindEQ__23_}   the following inequality holds 
\begin{align}\label{1}
V(t)\le V(0), \quad  \forall t\in \left[0,{\mathop{\lim }\limits_{i\to +\infty }} \left(t_{i} \right)\right).\end{align}
Equations \eqref{GrindEQ__3_}, \eqref{GrindEQ__4_},  \eqref{GrindEQ__7_} and  \eqref{GrindEQ__9_} with $i= 0$, imply that $\textcolor{black}{\rho(0,0)=\rho_s}$  or equivalently that 
\begin{align}\label{2} \textcolor{black}{\left|\ln \left(\frac{\rho (0,0)}{\rho _{s} } \right)\right|=0.}\end{align}
Applying \eqref{GrindEQ__19_}, for $t=0$  and using \eqref{GrindEQ__5_}, and \eqref{2}, we obtain: 
\begin{align} \label{3} 
V(0)&={\mathop{\sup }\limits_{\textcolor{black}{0< x\le 1}}} \left(\left|\ln \left(\frac{\rho_0 (x)}{\rho _{s} } \right)\right|\exp (-\sigma x)\right)\nonumber\\ &\leq {\mathop{\sup }\limits_{\textcolor{black}{0< x\le 1}}} \left(\left|\ln \left(\frac{\rho_0 (x)}{\rho _{s} } \right)\right|\right).
\end{align} 
On the other hand, definition \eqref{GrindEQ__19_} implies for all $t\in \left[0,{\mathop{\lim }\limits_{i\to +\infty }} \left(t_{i} \right)\right)$
\begin{align} \label{4} 
V(t)&={\mathop{\sup }\limits_{0\le x\le 1}} \left(\left|\ln \left(\frac{\rho (t,x)}{\rho _{s} } \right)\right|\exp (-\sigma x)\right)\nonumber\\ & \geq \exp (-\sigma ) {\mathop{\sup }\limits_{0\le x\le 1}} \left(\left|\ln \left(\frac{\rho (t,x)}{\rho _{s} } \right)\right|\right).
\end{align} 
Combining \eqref{1}, \eqref{3} and \eqref{4}, we get 
 the following estimate for all $t\in \left[0,{\mathop{\lim }\limits_{i\to +\infty }} \left(t_{i} \right)\right)$:
\begin{align} \label{GrindEQ__24_} 
{\mathop{\sup }\limits_{0\le x\le 1}} \left(\left|\ln \left(\frac{\rho (t,x)}{\rho _{s} } \right)\right|\right)\le \exp (\sigma ){\mathop{\sup }\limits_{0<x\le 1}} \left(\left|\ln \left(\frac{\rho _{0} (x)}{\rho _{s} } \right)\right|\right).
\end{align} 
 Definition \eqref{GrindEQ__12_} in conjunction with \eqref{GrindEQ__2_}, \eqref{GrindEQ__24_} and the facts that $v(t)=\lambda \left(W(t)\right)$, $\lambda \in C^{1} \left(\mathbb{R} _{+} ;(0,+\infty )\right)$ is a non-increasing function imply that $v(t)\ge c(\rho _{0} )$ for all $t\in \left[0,{\mathop{\lim }\limits_{i\to +\infty }} \left(t_{i} \right)\right)$. Estimate \eqref{GrindEQ__10_} for $t\in \left[0,{\mathop{\lim }\limits_{i\to +\infty }} \left(t_{i} \right)\right)$ is a direct consequence of estimate \eqref{GrindEQ__23_}, definition \eqref{GrindEQ__19_} and the fact that $v(t)\ge c(\rho _{0} )$ for all $t\in \left[0,{\mathop{\lim }\limits_{i\to +\infty }} \left(t_{i} \right)\right)$.

 The rest of the proof is devoted to the proof of \eqref{GrindEQ__11_} which also shows that ${\mathop{\lim }\limits_{i\to +\infty }} \left(t_{i} \right)=+\infty $.
 Define for all $t\in [t_{i} ,t_{i+1} )$ and $i\ge 0$:
\begin{align} \label{GrindEQ__25_} 
z(t)=\exp \left(\sigma \int _{t_{i} }^{t}\lambda \left(W(s)\right)ds \right)\ln \left(\frac{\rho (t,0)}{\rho _{s} } \right) 
\end{align} 
By virtue of \eqref{GrindEQ__25_}, \eqref{GrindEQ__3_}, \eqref{GrindEQ__4_} and \eqref{GrindEQ__9_}, we get for all $t\in [t_{i} ,t_{i+1} )$ and $i\ge 0$:
\begin{align} \label{GrindEQ__26_} 
z(t)=\exp \left(\sigma \int _{t_{i} }^{t}\lambda \left(W(s)\right)ds \right)\ln \left(\frac{\lambda \left(W(t_{i} )\right)}{\lambda \left(W(t)\right)} \right) 
\end{align} 

\noindent Equation \eqref{GrindEQ__26_} implies for all $i\ge 0$ and  for $t\in [t_{i} ,t_{i+1}) $ \textcolor{black}{almost everywhere}:
\begin{align} \label{GrindEQ__27_}
\dot{z}(t)&=-\lambda '\left(W(t)\right)\left(\rho (t,0)-\rho (t,1)\right)\nonumber\\ &\times \exp \left(\sigma \int _{t_{i} }^{t}\lambda \left(W(s)\right)ds \right)+\sigma \lambda \left(W(t)\right)z(t).             
\end{align} 
 Differential equation \eqref{GrindEQ__27_} is a direct consequence of the fact that for $t\in [t_{i} ,t_{i+1} )$
  $$W(t)=(t-t_{i} )u_{i} +\int _{0}^{a(t)}\rho (t_{i} ,x)dx,$$ 
 with $$a(t)=1-\int _{t_{i} }^{t}\lambda (W(s))ds.$$ The previous equations and the fact that $\rho [t_{i} ]$ is piecewise continuous show that the mapping $[t_{i} ,t_{i+1} )\ni t\to W(t)$ is continuous and piecewise continuously differentiable for $t\in [t_{i} ,t_{i+1} )$. As a consequence of \eqref{GrindEQ__26_}, we also obtain that the mapping $[t_{i} ,t_{i+1} )\ni t\to z(t)$ is continuous and piecewise continuously differentiable for $t\in [t_{i} ,t_{i+1} )$.

 It follows from the triangle inequality, \eqref{GrindEQ__25_}, \eqref{GrindEQ__27_}, the fact that $\lambda \in C^{1} \left(\mathbb{R} _{+} ;(0,+\infty )\right)$ is a non-increasing function and the fact that $\left|\lambda '(s)\right|\le K$ for all $s\ge 0$, that the following inequality holds for all $i\ge 0$ and $t\in [t_{i} ,t_{i+1} )$ almost everywhere:
\begin{align} \label{GrindEQ__28_} 
& \left|\dot{z}(t)\right| \le K\left|\rho (t,0)-\rho (t,1)\right|\exp \left(\sigma \int _{t_{i} }^{t}\lambda \left(W(s)\right)ds \right)\nonumber \\&+\sigma \lambda \left(0\right)\left|z(t)\right| \\ &\le K\left|\rho (t,0)-\rho _{s} \right|\exp \left(\sigma \int _{t_{i} }^{t}\lambda \left(W(s)\right)ds \right)\nonumber\\ &+K\left|\rho (t,1)-\rho _{s} \right|\exp \left(\sigma \lambda (0)(t-t_{i} )\right)+\sigma \lambda \left(0\right)\left|z(t)\right|
\end{align} 
Formulas \eqref{GrindEQ__18_}, \eqref{GrindEQ__19_} and the fact that $\left|\exp (x)-1\right|\le \exp (\left|x\right|)-1$ for all $x\in \mathbb{R} $, leads to the following estimates 
\begin{align}\left|\rho (t,1)-\rho _{s} \right|&\le \left\| \rho [t_{i} ]-\rho _{s} \right\| _{\infty },\nonumber\\
\left\| \rho [t_{i} ]-\rho _{s} \right\| _{\infty }& \le \rho _{s} \left(\exp \left(\exp (\sigma )V(t_{i} )\right)-1\right).\nonumber\end{align}
 Moreover, the fact that $\left|\exp (x)-1\right|\le \left|x\right|\exp (\left|x\right|)$ for all $x\in \mathbb{R} $, shows that $$\left|\rho (t,0)-\rho _{s} \right|\le \rho _{s} \exp \left(\left|\ln \left(\frac{\rho (t,0)}{\rho _{s} } \right)\right|\right)\left|\ln \left(\frac{\rho (t,0)}{\rho _{s} } \right)\right|.$$ Using the previous inequalities in conjunction with \eqref{GrindEQ__19_}, \eqref{GrindEQ__22_}, \eqref{GrindEQ__28_} and \eqref{GrindEQ__25_}, we get the following inequality for all $i\ge 0$ and $t\in [t_{i} ,t_{i+1} )$ almost everywhere:
\begin{align} \label{GrindEQ__29_} 
&\left|\dot{z}(t)\right|\le \left(K\rho _{s} \exp \left(V(t_{i} )\right)+\sigma \lambda \left(0\right)\right)\left|z(t)\right|\nonumber\\&+K\rho _{s} \left(\exp \left(\exp (\sigma )V(t_{i} )\right)-1\right)\exp \left(\sigma \lambda (0)(t-t_{i} )\right). 
\end{align} 
Using \eqref{GrindEQ__29_}, the facts that $z(t_{i} )=0$, $t_{i+1} -t_{i} \le \frac{1}{\lambda (0)} $ (a direct consequence of (8)),  and the Gronwall-Bellman lemma, we get for all $i\ge 0$ and $t\in [t_{i} ,t_{i+1} )$: 
\begin{align} \label{GrindEQ__30_} 
\left|z(t)\right|&\le \left(\exp \left(\exp (\sigma )V(t_{i} )\right)-1\right)\nonumber\\ &\times \exp \left(\sigma \right)\frac{\exp \left(K\rho _{s} \exp \left(V(t_{i} )\right)(t-t_{i} )\right)-1}{\exp \left(V(t_{i} )\right)}. 
\end{align} 
 It follows from \eqref{GrindEQ__30_} that the inequality $\left|z(t)\right|\le V(t_{i} )$ holds for all $t\in [t_{i} ,t_{i+1} )$, provided that
\begin{align} \label{GrindEQ__31_} 
t_{i+1} -t_{i} \le \tilde{T}\left(V(t_{i} )\right), 
\end{align} 
where 
\begin{equation} \label{GrindEQ__32_} 
\tilde{T}(s):=\left\{\begin{array}{c} {\frac{1}{K\rho _{s} } \exp \left(-s\right)\ln \left(1+\frac{s\exp \left(s-\sigma \right)}{\exp \left(\exp (\sigma )s\right)-1} \right)\, if\, s>0} \\ {\frac{1}{K\rho _{s} } \ln \left(1+\exp \left(-2\sigma \right)\right)\quad if\quad s=0} \end{array}\right.
\end{equation} 
with  $\tilde T(0)=\mathop{\lim } \limits_{s \to 0^+ } \tilde T \left(s \right).$ 

Using \eqref{GrindEQ__8_a}, \eqref{GrindEQ__8_b}, \eqref{GrindEQ__25_} and \eqref{GrindEQ__32_}, we conclude that 
\begin{align} \label{GrindEQ__33_} 
t_{i+1} -t_{i} \ge \min \left(\, \frac{1}{\lambda (0)} \, ,\, \tilde{T}(V(t_{i} ))\, \right). 
\end{align} 
 Continuity and positivity of the function $\tilde{T}(s)$ defined by \eqref{GrindEQ__32_}, implies that there exists a non-increasing function $T:\mathbb{R} _{+} \to (0,+\infty )$ such that $$\min \left(\, \frac{1}{\lambda (0)} \, ,\, \tilde{T}(s)\, \right)\ge T(s),$$ for all $s\ge 0$. Inequality \eqref{GrindEQ__11_} is a consequence of the previous inequality, \eqref{GrindEQ__33_}, \eqref{GrindEQ__23_} and definition \eqref{GrindEQ__19_}. 
\noindent The proof is complete.      $\triangleleft $ \\
\begin{rem}
 Estimate \eqref{GrindEQ__10_} is a stability estimate in a special state norm. Due to the positivity of the state, the logarithmic norm of the state $\rho $ appears, i.e., we have $\left|\ln \left(\frac{\rho (t,x)}{\rho _{s} } \right)\right|$ instead of the usual $\left|\rho (t,x)-\rho _{s} \right|$ that appears in many stability estimates for linear PDEs. The logarithmic norm is a manifestation of the nonlinearity of system \eqref{GrindEQ__1_}, \eqref{GrindEQ__2_}, \eqref{GrindEQ__3_} and the fact that the state space is not a linear space but rather a positive cone: the state space for system \eqref{GrindEQ__1_}, \eqref{GrindEQ__2_}, \eqref{GrindEQ__3_} is the set $X:=\left\{\, \rho \in PC^{1} \left([0,1]\right)\, :\, {\mathop{\inf }\limits_{x\in (0,1]}} \left(\rho (x)\right)>0\, \right\}$. The use of the logarithmic norm of the state is common in systems with positivity constraints (see \textcolor{black}{\cite{s33,s33x}}). 
\end{rem}
\begin{rem}
 Estimate \eqref{GrindEQ__11_} guarantees that no Zeno behavior can appear for the closed-loop system \eqref{GrindEQ__1_}, \eqref{GrindEQ__2_}, \eqref{GrindEQ__3_}, \eqref{GrindEQ__4_}, \eqref{GrindEQ__7_}, \eqref{GrindEQ__8_a}, \eqref{GrindEQ__8_b}, \eqref{GrindEQ__9_}. 
 The proof of Theorem 2 provides an estimate for the function $T:\mathbb{R} _{+} \to (0,+\infty )$. 
\end{rem}
\subsection{Robustness With Respect to the Event Sequence}
\noindent Let $\rho _{0} \in PC^{1} \left([0,1]\right)$ with ${\mathop{\inf }\limits_{x\in (0,1]}} \left(\rho _{0} (x)\right)>0$ be given and define:
\begin{align} \label{GrindEQ__13_} 
G(\rho _{0} ):=\min \left(\, \frac{1}{\lambda (0)} \, ,\, r\, \right),
\end{align} 
where 
\begin{align} \label{GrindEQ__14_} 
r:=\inf &\Bigg\{\, \tau >0\, : \left|\ln \left(\frac{\rho (\tau ,0)}{\rho _{s} } \right)\right|\nonumber\\ &>\exp \left(-\sigma \int _{0}^{\tau }\lambda (W(s))ds \right)\nonumber \\&\times {\mathop{\sup }\limits_{0<x\le 1}} \left(\left|\ln \left(\frac{\rho _{0} (x)}{\rho _{s} } \right)\right|\exp (-\sigma x)\right)\Bigg\}, 
\end{align} 
 where $\rho [t]\in PC^{1} \left([0,1]\right)$ and ${\mathop{\inf }\limits_{x\in [0,1]}} \left(\rho (t,x)\right)>0$ for all $t\ge 0$, is the solution of \eqref{GrindEQ__1_}, \eqref{GrindEQ__2_}, \eqref{GrindEQ__3_}, \eqref{GrindEQ__5_} with $u(t)\equiv \rho _{s} \lambda \left(\int _{0}^{1}\rho _{0} (x)dx \right)$. Definitions \eqref{GrindEQ__13_}, \eqref{GrindEQ__14_} imply that the event-triggered control \eqref{GrindEQ__7_}, \eqref{GrindEQ__8_a}, \eqref{GrindEQ__8_b}, \eqref{GrindEQ__9_} satisfies the following relation for all $i\ge 0$:
\begin{align} \label{GrindEQ__15_} 
t_{i+1} =t_{i} +G\left(\rho [t_{i} ]\right) 
\end{align} 
Theorem \ref{T2} (and particularly inequality (11)) shows that for every initial condition $\rho _{0} \in PC^{1} \left([0,1]\right)$ with ${\mathop{\inf }\limits_{x\in (0,1]}} \left(\rho _{0} (x)\right)>0$, the sequence of events $\left\{\, t_{i} \, :\, i=0,1,2,...\, \right\}$ with $t_{0} =0$ is a diverging sequence, i.e., ${\mathop{\lim }\limits_{i\to +\infty }} \left(t_{i} \right)=+\infty $. However, there is an infinite number of diverging increasing sequences $\left\{\, t_{i} \, :\, i=0,1,2,...\, \right\}$ with $t_{0} =0$ for which $t_{i+1} -t_{i} \le G\left(\rho [t_{i} ]\right)$ for all $i\ge 0$. For these sequences, the controller acts (through (9)) before an event occurs. The following result extends the result of Theorem 2 and guarantees robustness with respect to the event sequence.
\begin{theorem}\label{T3}Suppose that there exists a constant $K>0$ such that $\left|\lambda '(s)\right|\le K$ for all $s\ge 0$. Let $\sigma >0$ be a given parameter. Then for every $\rho _{0} \in PC^{1} \left([0,1]\right)$ with ${\mathop{\inf }\limits_{x\in (0,1]}} \left(\rho _{0} (x)\right)>0$ and for every increasing sequence of times $\left\{\, t_{i} \, :\, i=0,1,2,...\, \right\}$ with $t_{0} =0$, ${\mathop{\lim }\limits_{i\to +\infty }} \left(t_{i} \right)=+\infty $ that satisfies
\begin{equation} \label{GrindEQ__16_}
t_{i+1} -t_{i} \le G\left(\rho [t_{i} ]\right), \ for \ all \ i\ge 0,                                                
\end{equation} 
the solution $\rho :\mathbb{R} _{+} \times [0,1]\to (0,+\infty )$ with $\rho [t]\in PC^{1} \left([0,1]\right)$ and ${\mathop{\inf }\limits_{x\in [0,1]}} \left(\rho (t,x)\right)>0$ for all $t\ge 0$ of the initial-boundary value problem \eqref{GrindEQ__1_}, \eqref{GrindEQ__2_}, \eqref{GrindEQ__3_}, \eqref{GrindEQ__4_}, \eqref{GrindEQ__5_}, \eqref{GrindEQ__9_} satisfies estimate \eqref{GrindEQ__10_}. \end{theorem}

\textbf{Proof of Theorem \ref{T3}:} The proof is essentially the same as the first part of the proof of Theorem 2. The only difference is to notice that the event-trigger \eqref{GrindEQ__8_a}, \eqref{GrindEQ__8_b}, (as well as continuity of $v(t)=\lambda \left(W(t)\right)$, which implies continuity of $\rho (t,0)=\frac{u_{i} }{v(t)} $ for $t\in \left[t_{i} ,\min (t_{i+1} ,t_{\max } )\right)$) gives for all $i\ge 0$ with $t_{i} <t_{\max } $ when combined with definitions \eqref{GrindEQ__19_}, \eqref{GrindEQ__13_}, \eqref{GrindEQ__14_}:
\begin{align} \label{GrindEQ__34_}
\left|\ln \left(\frac{\rho (\tau ,0)}{\rho _{s} } \right)\right|\le \exp \left(-\sigma \int _{t_{i} }^{\tau }\lambda (W(s))ds \right)V(t_{i} ), \end{align} 
for all $\tau \in \left[t_{i} ,\min \left(t_{\max } ,t_{i} +G\left(\rho [t_{i} ]\right)\right)\right).$        
 Inequality \eqref{GrindEQ__34_} replaces inequality \eqref{GrindEQ__21_}. Using \eqref{GrindEQ__34_} we end up with inequality \eqref{GrindEQ__24_} exactly as in the proof of Theorem 2. The rest of the proof of Theorem 2 is not needed because it is not needed to show inequality \eqref{GrindEQ__11_} and  it is assumed that ${\mathop{\lim }\limits_{i\to +\infty }} \left(t_{i} \right)=+\infty $ . The proof is complete.       $\triangleleft $ \\

\subsection{ Sampled-Data Stabilization with Robustness With Respect to the Sampling Schedule}

\noindent Theorem \ref{T3} is important because it shows that the controller \eqref{GrindEQ__9_} can be implemented in various ways. For example, we can implement the controller \eqref{GrindEQ__9_} in a sample-and-hold fashion for an appropriate sampling period. This is shown by the following result. 

\begin{theorem}\label{T4}Suppose that there exists a constant $K>0$ such that $\left|\lambda '(s)\right|\le K$ for all $s\ge 0$. Let $\sigma >0$ be a given parameter. Then, for every $R>0$ there exists $\tau >0$ such that for every $\rho _{0} \in PC^{1} \left([0,1]\right)$ with ${\mathop{\inf }\limits_{x\in (0,1]}} \left(\rho _{0} (x)\right)>0$, ${\mathop{\sup }\limits_{0<x\le 1}} \left(\left|\ln \left(\frac{\rho _{0} (x)}{\rho _{s} } \right)\right|\right)\le R$ and for every increasing sequence of times $\left\{\, t_{i} \, :\, i=0,1,2,...\, \right\}$ with $t_{0} =0$, ${\mathop{\lim }\limits_{i\to +\infty }} \left(t_{i} \right)=+\infty $ that satisfies
\begin{align} \label{GrindEQ__17_}
t_{i+1} -t_{i} \le \tau , \ \textrm{for} \ \textrm{all} \ i\ge 0,                                                   
\end{align} 
the solution $\rho :\mathbb{R} _{+} \times [0,1]\to (0,+\infty )$ with $\rho [t]\in PC^{1} \left([0,1]\right)$ and ${\mathop{\inf }\limits_{x\in [0,1]}} \left(\rho (t,x)\right)>0$ for all $t\ge 0$ of the initial-boundary value problem \eqref{GrindEQ__1_}, \eqref{GrindEQ__2_}, \eqref{GrindEQ__3_}, \eqref{GrindEQ__4_}, \eqref{GrindEQ__5_}, \eqref{GrindEQ__9_} satisfies estimate \eqref{GrindEQ__10_}. \end{theorem}
 The sample-and-hold implementation of the controller \eqref{GrindEQ__9_} does not require continuous measurement of the state. On the other hand, the time $\tau >0$ is (in general) much smaller than $G\left(\rho [t_{i} ]\right)$, which implies that the control action must be updated much more frequently in the sampled-data case than in the event-triggered case. 
 
\textbf{\textcolor{black}{Proof of Theorem \ref{T4}:}} The proof of Theorem \ref{T2} and definitions \eqref{GrindEQ__13_}, \eqref{GrindEQ__14_} actually show that there exists a non-increasing function $T:\mathbb{R} _{+} \to (0,+\infty )$ such that for every $R>0$ and for every $\rho \in PC^{1} \left([0,1]\right)$ with ${\mathop{\inf }\limits_{x\in (0,1]}} \left(\rho (x)\right)>0$ and ${\mathop{\sup }\limits_{0<x\le 1}} \left(\left|\ln \left(\frac{\rho (x)}{\rho _{s} } \right)\right|\right)\le R$, the following inequality holds:
\begin{align} \label{GrindEQ__35_} 
G(\rho )\ge T\left(R\right).
\end{align} 
Setting $\tau :=T\left(R\right)$ and repeating the proof of Theorem 2 with \eqref{GrindEQ__34_} replacing \eqref{GrindEQ__21_}, we are in a position to show that inequality \eqref{GrindEQ__24_} holds for every increasing sequence of times $\left\{\, t_{i} \, :\, i=0,1,2,...\, \right\}$ with $t_{0} =0$, ${\mathop{\lim }\limits_{i\to +\infty }} \left(t_{i} \right)=+\infty $ that satisfies \eqref{GrindEQ__17_}. The proof is complete.       $\triangleleft $ 
\section{Simulation Results}\label{s5}
\subsection{Event-triggered control}

\begin{figure}[htbp]
\centering
\includegraphics[width=0.5\textwidth]{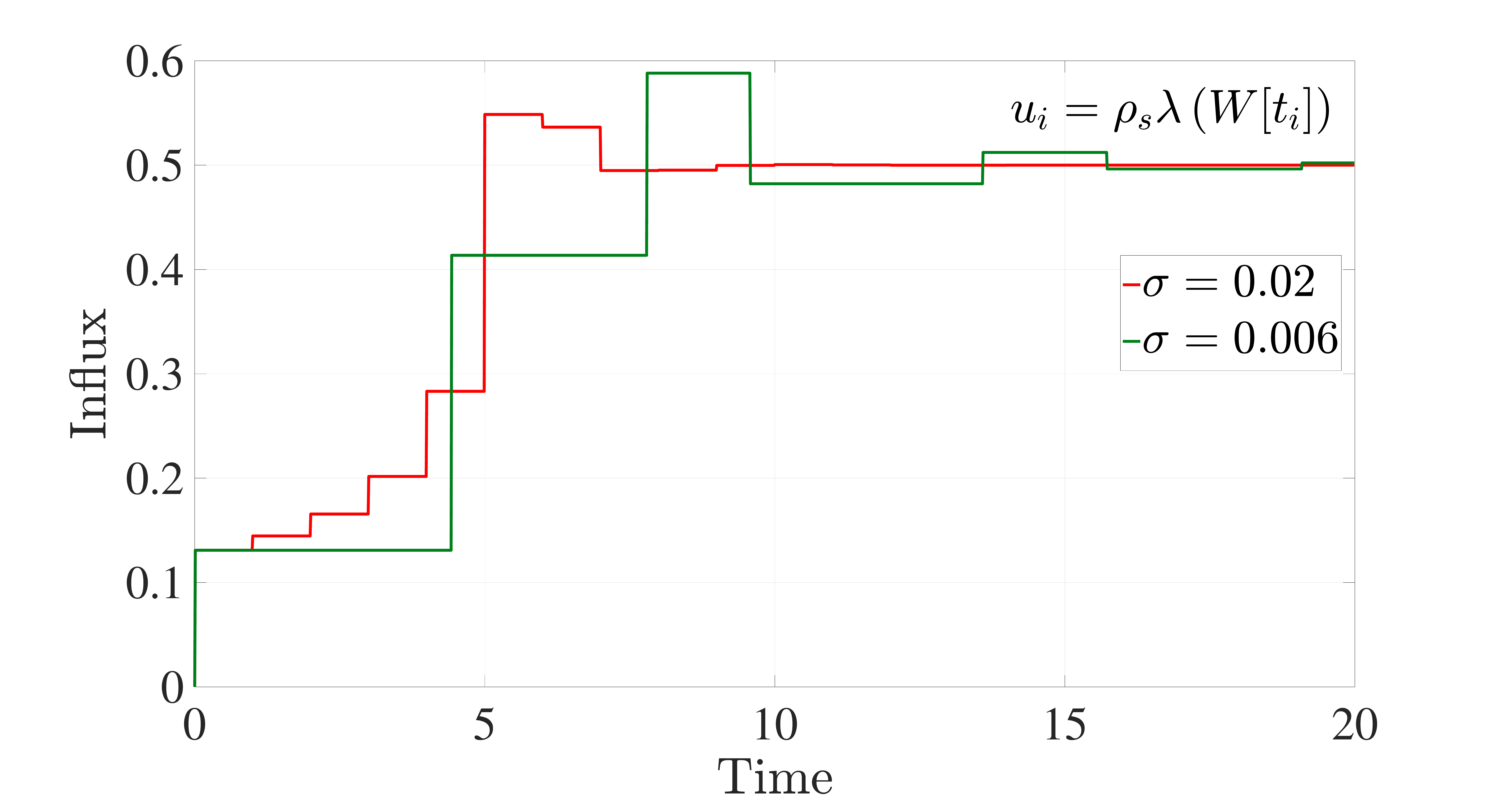}
\caption{Infux dynamics: event-triggered control action at the boundary $x=0$.}\label{system_figure0}
\end{figure}

\begin{figure}[htbp]
\centering
\includegraphics[width=0.5\textwidth]{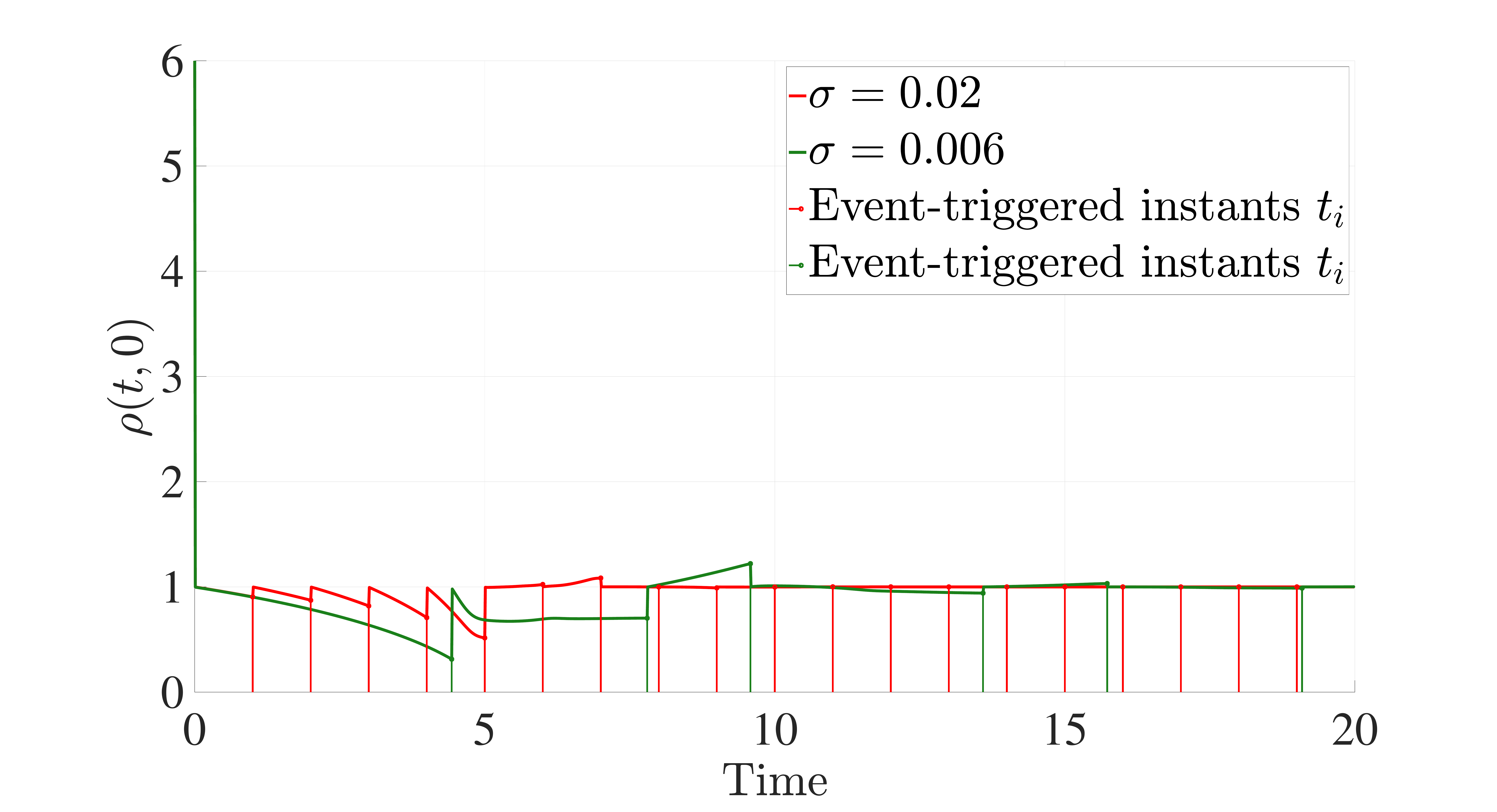}
\caption{Density evolution in time at the controlled boundary and  event-trigger instants  $t_i$.}\label{system_figure1}
\end{figure}

\begin{figure}[htbp]
\centering
\includegraphics[width=0.5\textwidth]{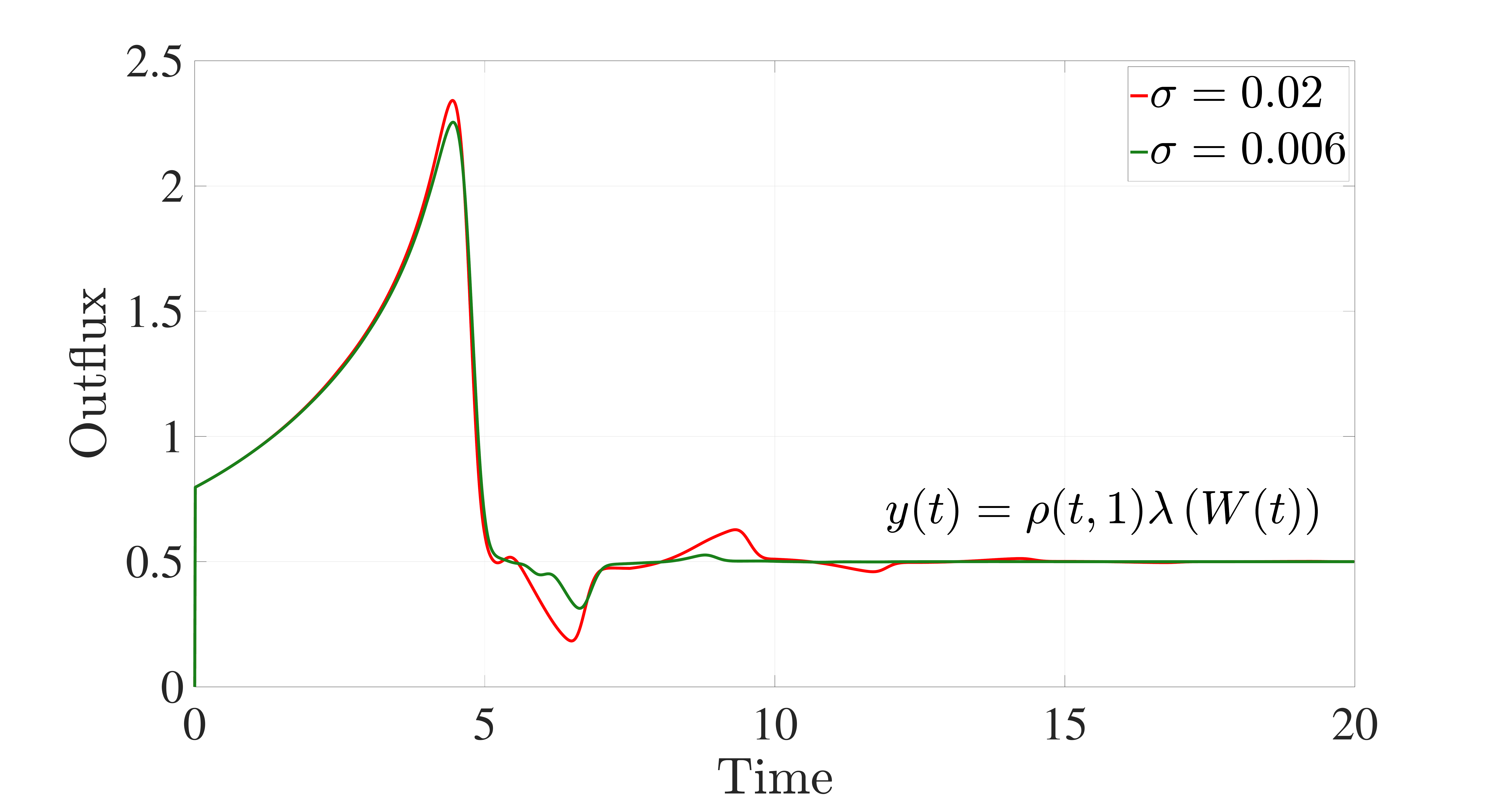}
\caption{Outflux dynamics at the uncontrolled boundary $x=1$.}\label{system_figure2}
\end{figure}

\begin{figure}[htbp]
\centering
\includegraphics[width=0.5\textwidth]{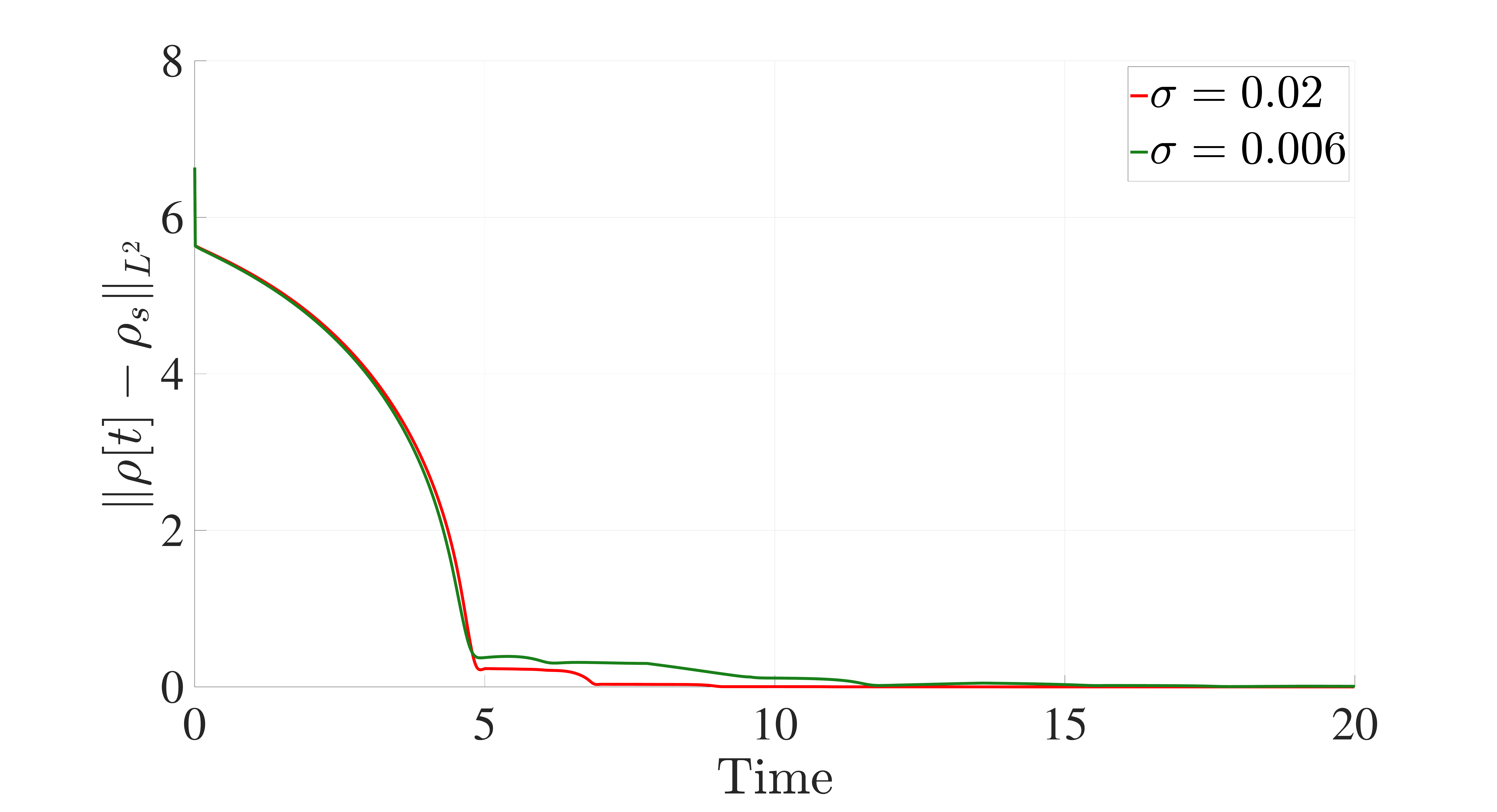}
\caption{$L^2$ norm of the distributed density deviation for different values of $\sigma$.}\label{system_figure3}
\end{figure}

\begin{figure}[htbp]
\centering
\includegraphics[width=0.45\textwidth]{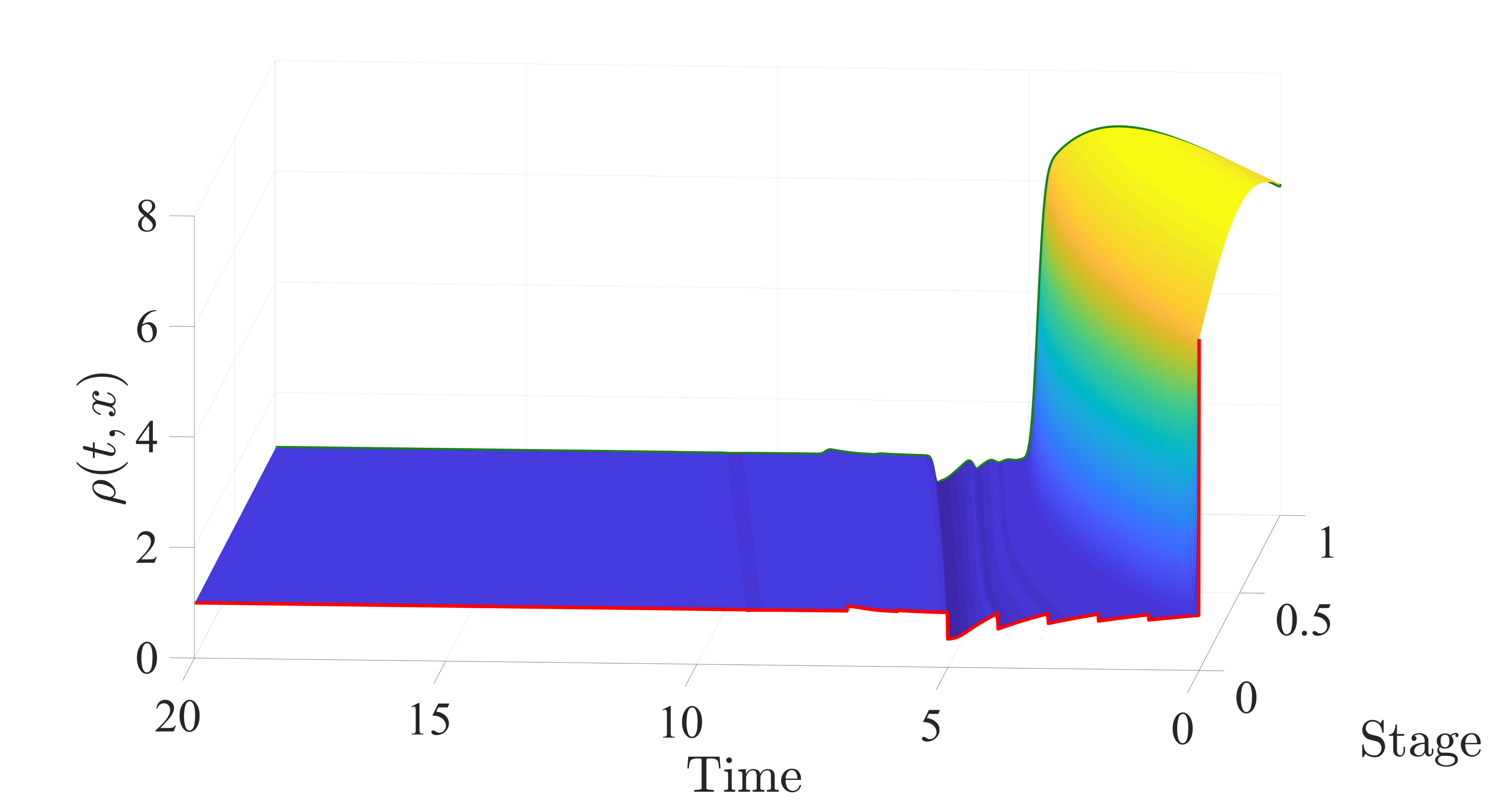}
\caption{The dynamics of the distributed density function for $\sigma=0.02$.}\label{system_figure4}
\end{figure}

\begin{figure}[htbp]
\centering
\includegraphics[width=0.45\textwidth]{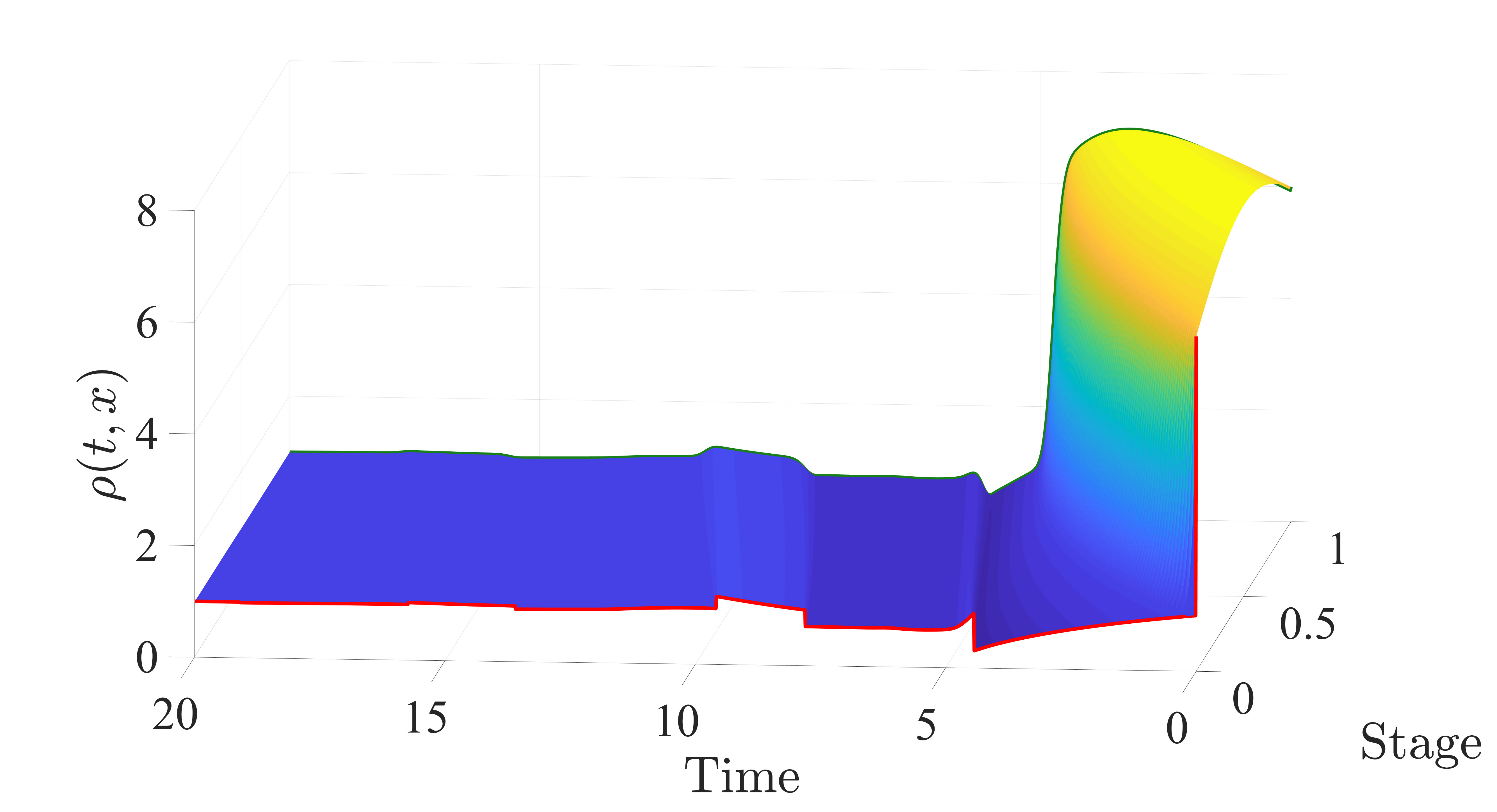}
\caption{The dynamics of the distributed density function  for $\sigma=0.006$.}\label{system_figure5}
\end{figure}

We simulate closed-loop system consisting of  \eqref{GrindEQ__1_}, 
\eqref{GrindEQ__2_}, \eqref{GrindEQ__3_}, \eqref{GrindEQ__5_} together with the even triggered controller \eqref{GrindEQ__7_}--\eqref{GrindEQ__9_}. The  nonlocal propagation speed of PDE \eqref{GrindEQ__1_} is defined as $$\lambda(W)=\frac{1}{1+W},$$ where $W$ is given by \eqref{GrindEQ__2_}. The initial condition is set to $\rho_0(x)=6+\sin(\pi x)$ and the equilibrium density is defined as $\rho_s=1$. The event generator is computed for two values of $\sigma$, namely, $\sigma=0.02$ and $\sigma=0.006$. Here, the number of updating times of the control signal is an increasing function of the parameter $\sigma$.  Figure \ref{system_figure0} shows the evolution of the influx of parts that is the boundary control action while Figure \ref{system_figure1} reflects the dynamics of the  density at the input with the event-triggered instants. From Figure \ref{system_figure1} and Figure \ref{system_figure2}, it can be viewed that both the input and the output density are stabilized at the desired uniform equilibrium $\rho_s=1$ for both values of $\sigma$. However, \textcolor{black}{the greater} is the value of $\sigma$, \textcolor{black}{the faster} is the convergence rate due to the increasing number of execution of the control task \eqref{GrindEQ__9_}. Moreover, Figure \ref{system_figure3} which represents the $L^2$-norm of   the deviation of the distributed state $\rho(t,x)$ with respect to the uniform equilibrium $\rho_s$ shows that $\rho(t,x)$ converges to the equilibrium in a $L^2$ sense. This statement is confirmed by Figure \ref{system_figure4} and Figure \ref{system_figure5} which shows the evolution in time of the distributed density $\rho(t,x)$ for both values of $\sigma$.

\begin{figure}[htbp]
\centering
\includegraphics[width=0.5\textwidth]{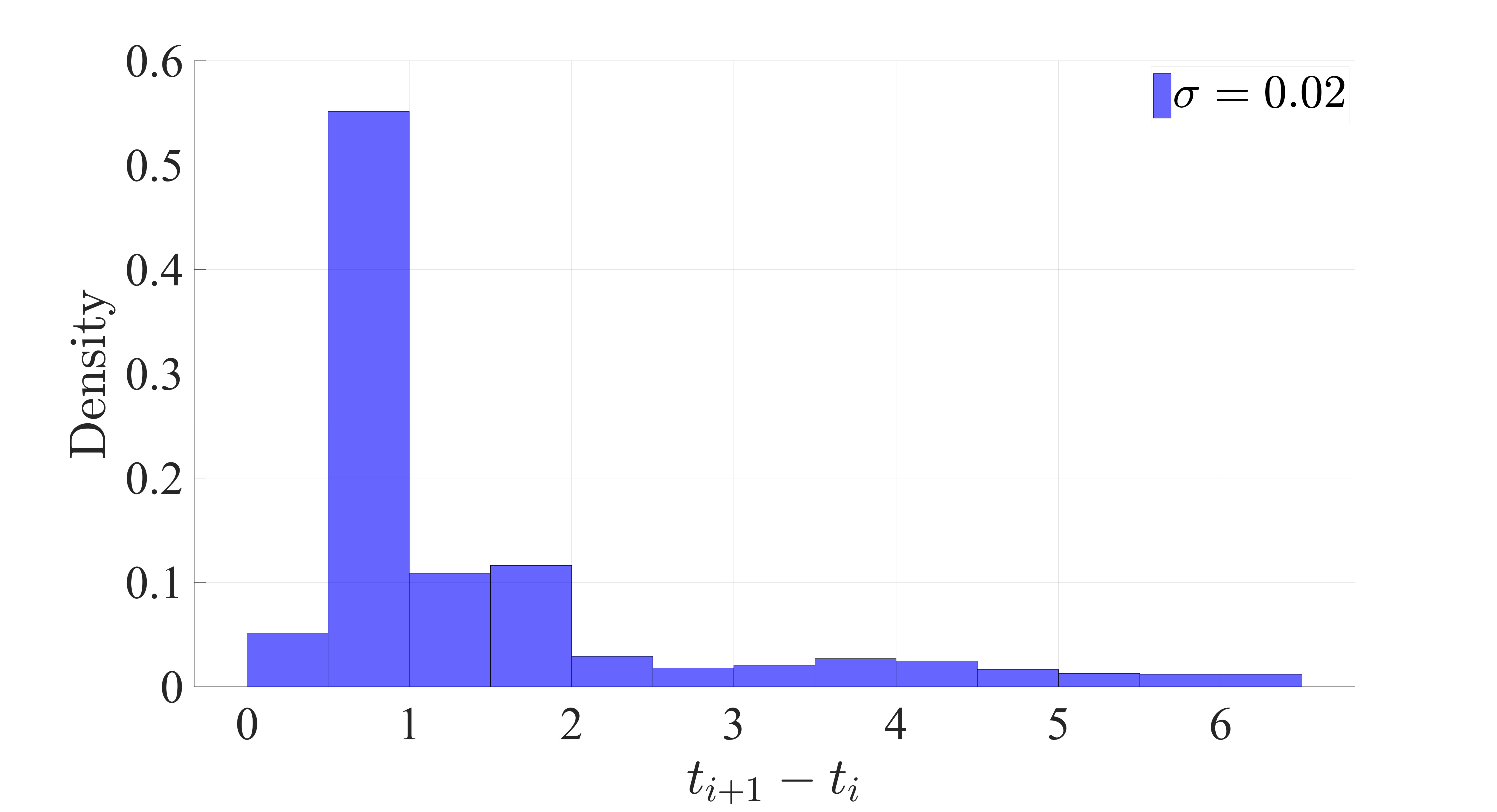}
\caption{The inter-execution time  $\sigma=0.02$: statistics with $200$ initial conditions.}\label{system_figure6}
\end{figure}

\begin{figure}[htbp]
\centering
\includegraphics[width=0.5\textwidth]{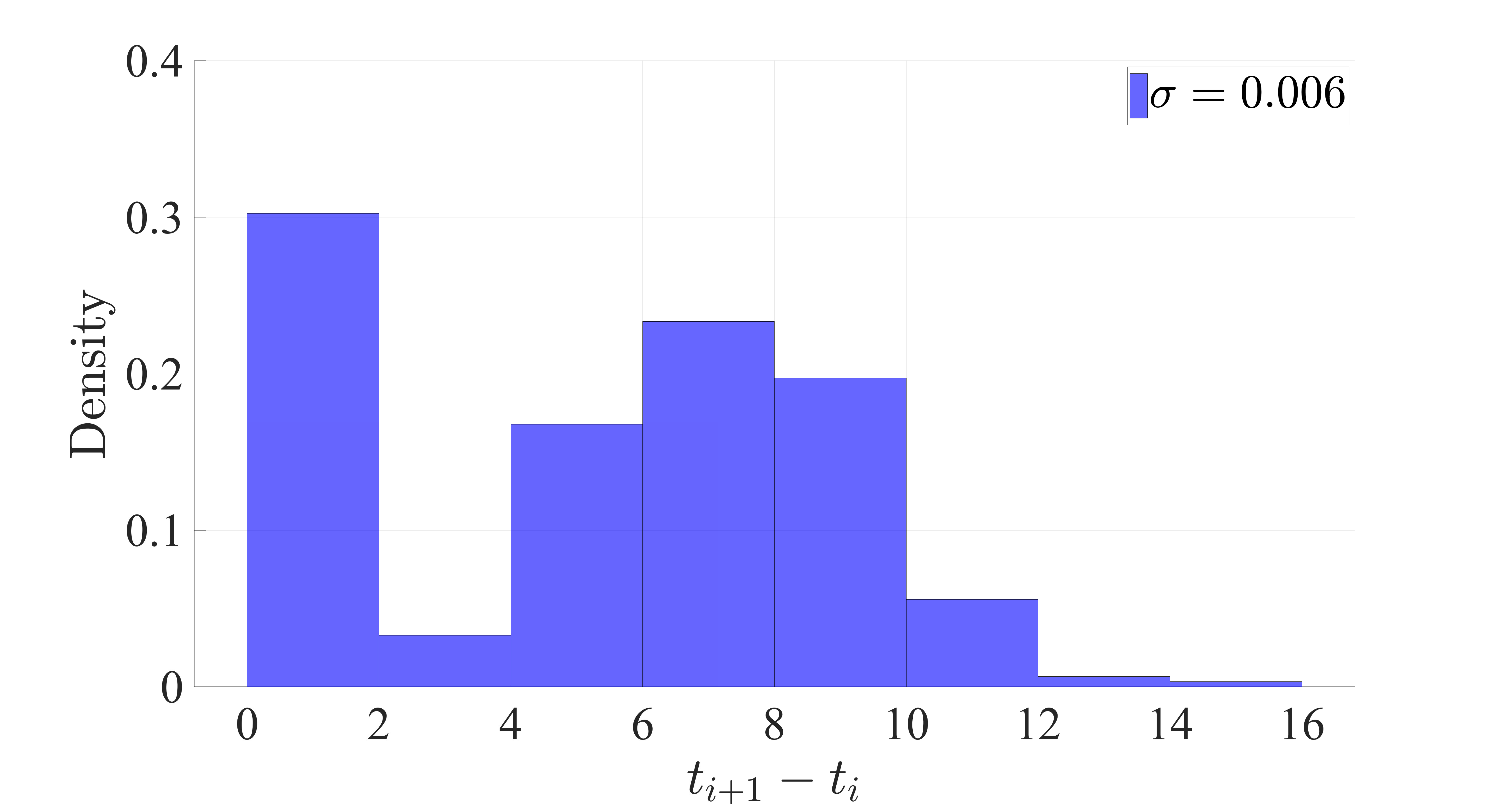}
\caption{The inter-execution time  $\sigma=0.006$: statistics on a $200$ initial conditions.}\label{system_figure7}
\end{figure}

 For a set of initial conditions defined as  $\rho_0(x)=6+\sin(\pi x)+ l\times x^4,\, l=1,\dots, 100$, Figure \ref{system_figure6} and  Figure \ref{system_figure7} represent the statistics on the inter-execution times under the event-triggered policy \eqref{GrindEQ__8_a}, \eqref{GrindEQ__8_b} for $\sigma=0.02$ and $\sigma=0.006$, respectively. \textcolor{black}{For a fast sampling, $\sigma=0.02$, the inter-execution time   $\tau=t_{i+1}-t_i$ belongs predominantly in the interval  $ [0.6,\ 1]$ while for slow sampling $\tau \in [0.6,\ 2]$ prevails. }These values can be used as indicative  (\textcolor{black}{but possibly conservative}) choices of the sampling periods  applying a sampled data control approach knowing the robustness of the event-based control with respect to the triggering policy. 

\subsection{Sampled data simulation results}
\begin{figure}[htbp]
\centering
\includegraphics[width=0.5\textwidth]{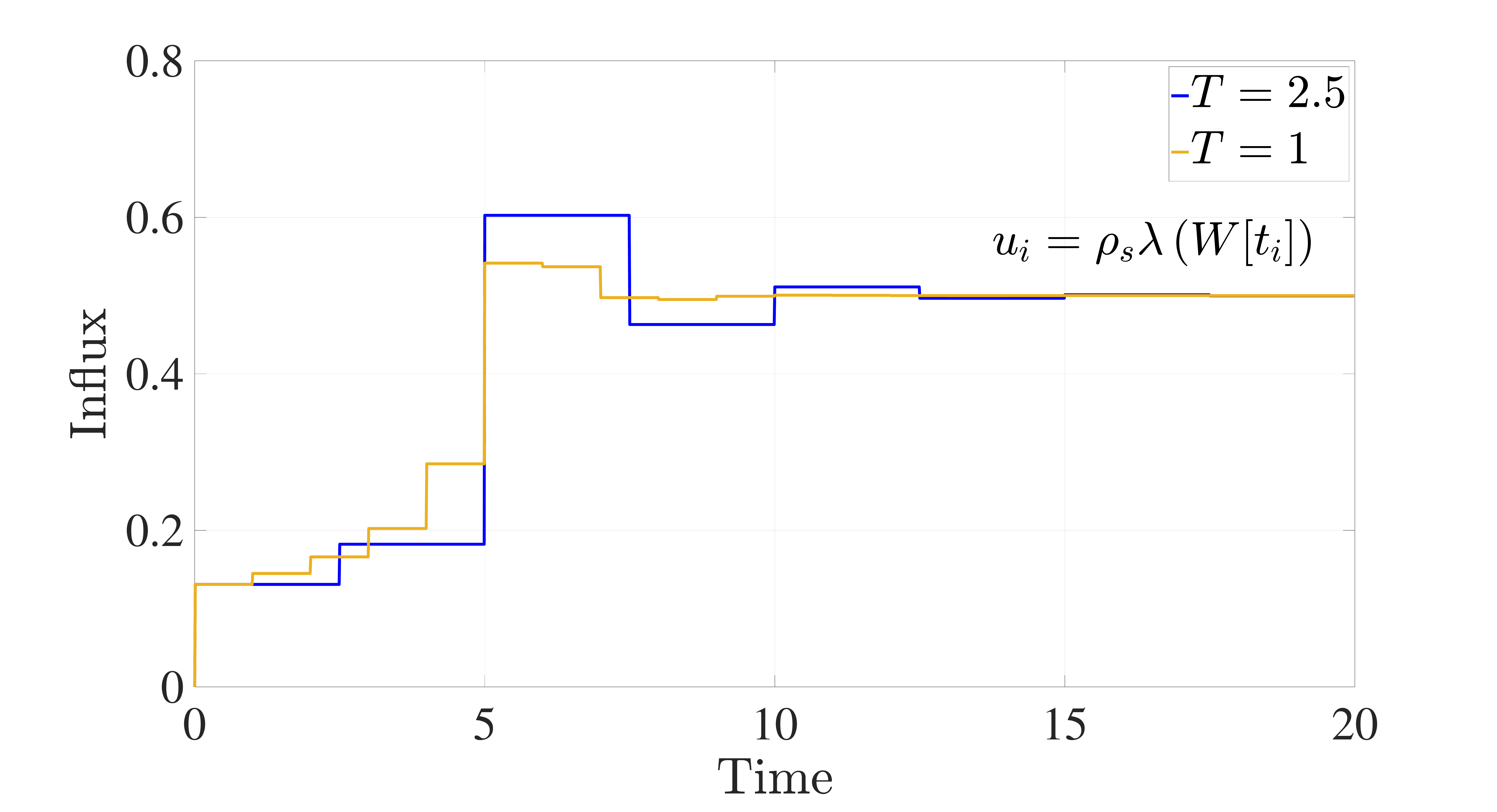}
\caption{Infux dynamics with sampled data boundary control action.}\label{system_figure-sample0}
\end{figure}

\begin{figure}[htbp]
\centering
\includegraphics[width=0.5\textwidth]{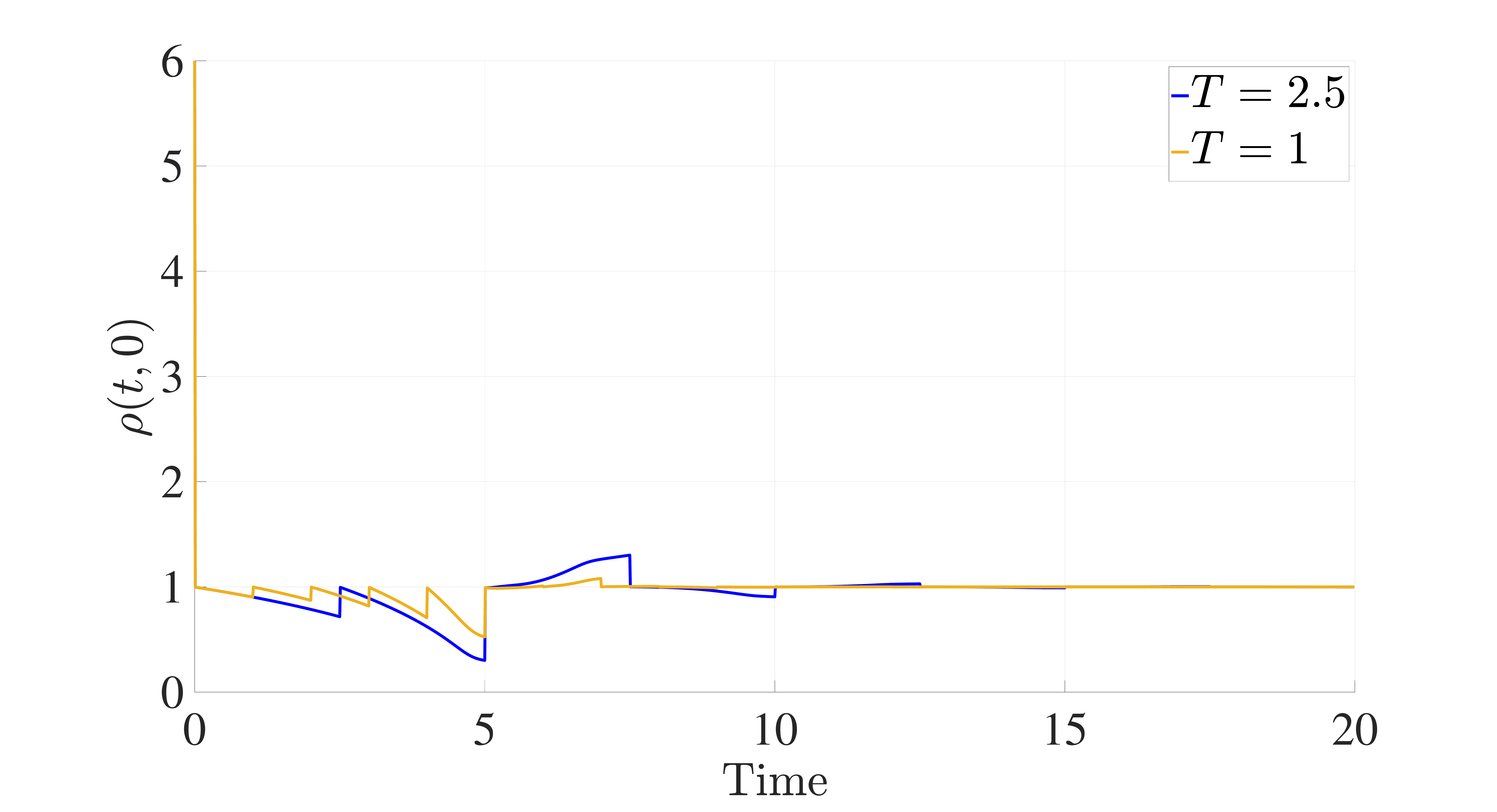}
\caption{The distributed density evolution in time at the controlled boundary ($x=0$) with sampled data control.}\label{system_figure-sample1} 
\end{figure}

\begin{figure}[htbp]
\centering
\includegraphics[width=0.5\textwidth]{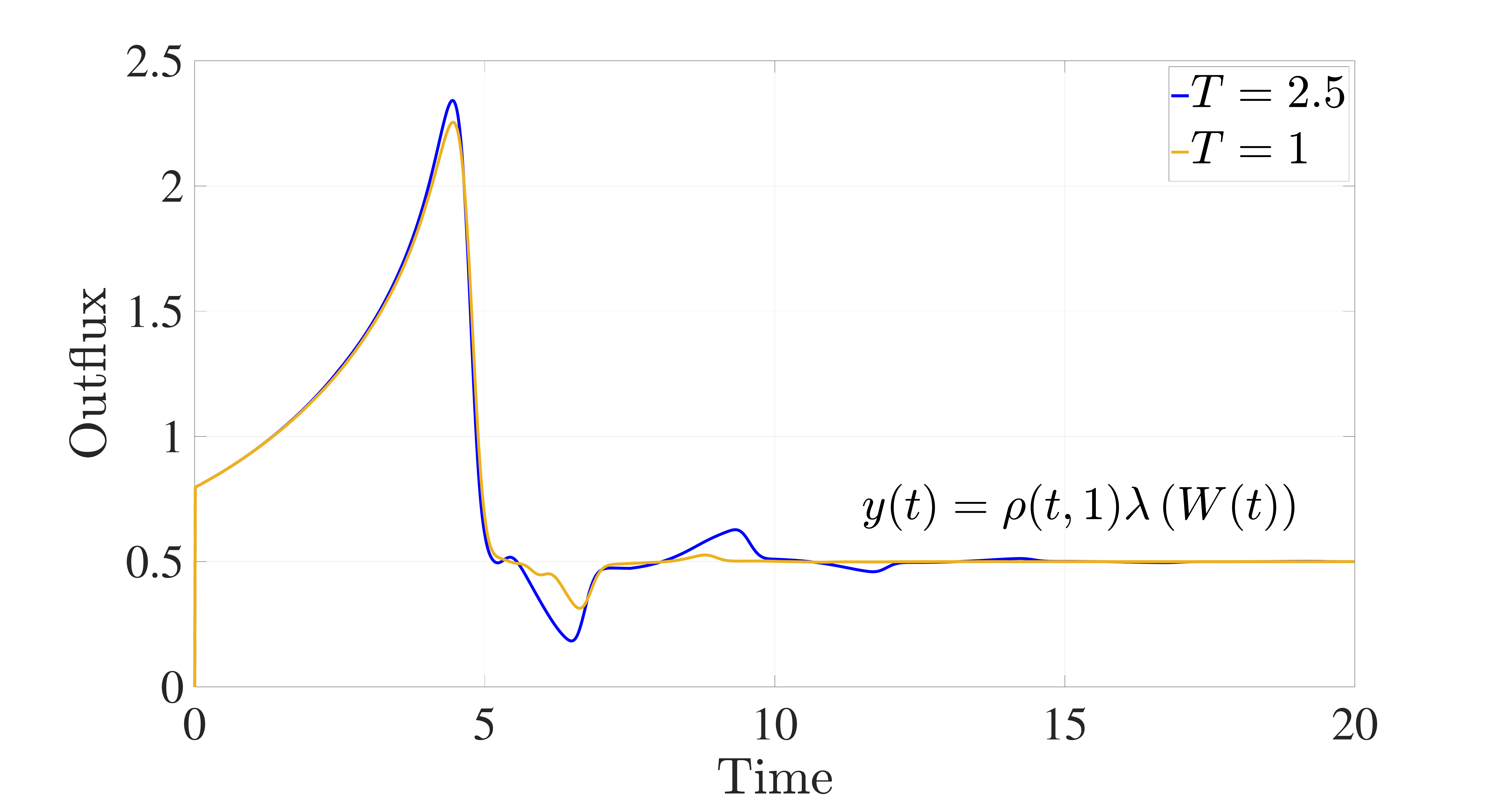}
\caption{ Outflux dynamics at the uncontrolled boundary $x=1$ with sampled data control.}\label{system_figure-sample2} 
\end{figure}

\begin{figure}[htbp]
\centering
\includegraphics[width=0.5\textwidth]{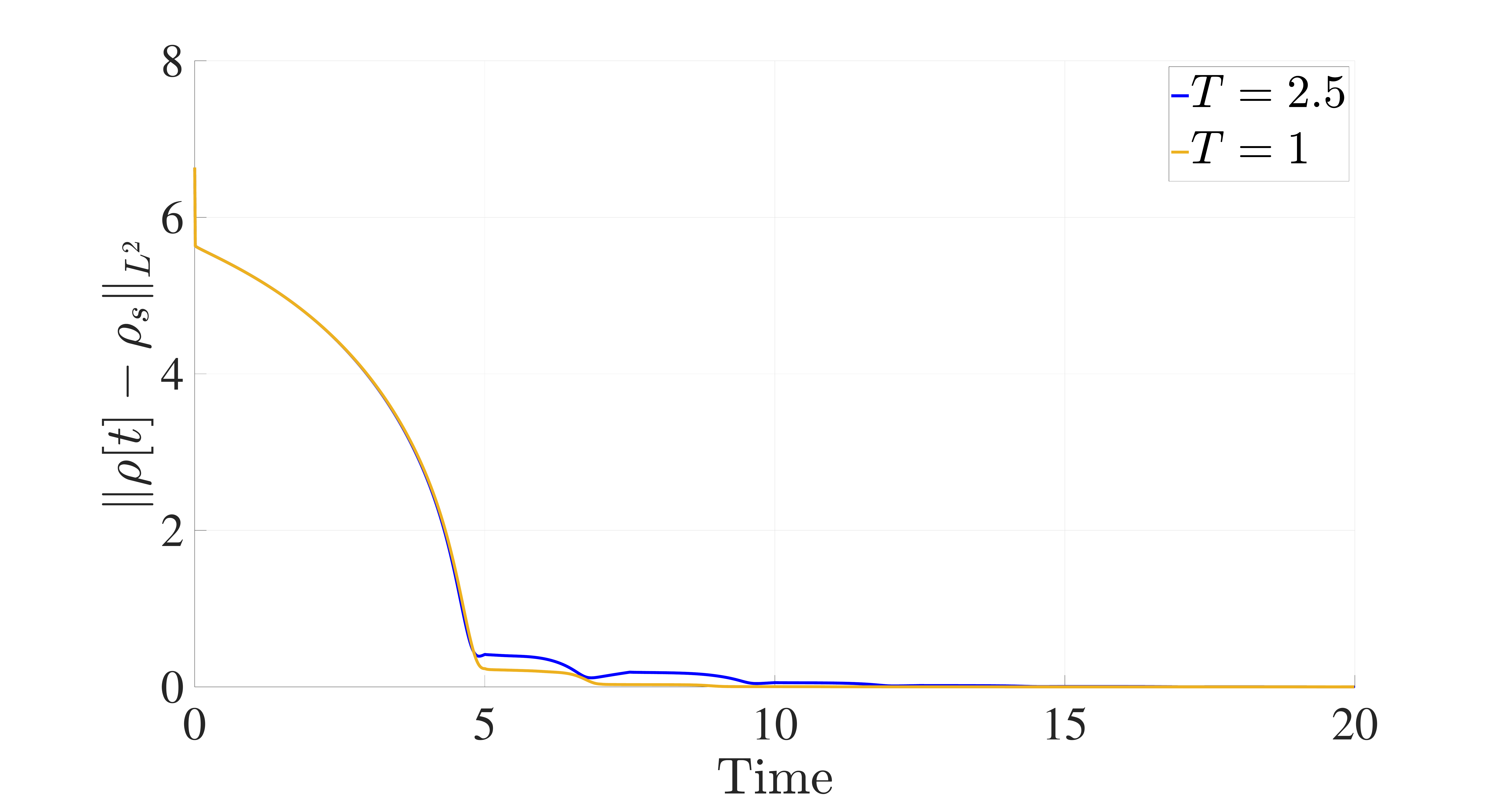}
\caption{$L^2$ norm of the distributed density deviation for the two sampling periods.}\label{system_figure-sample3} 
\end{figure}

\begin{figure}[htbp]
\centering
\includegraphics[width=0.45\textwidth]{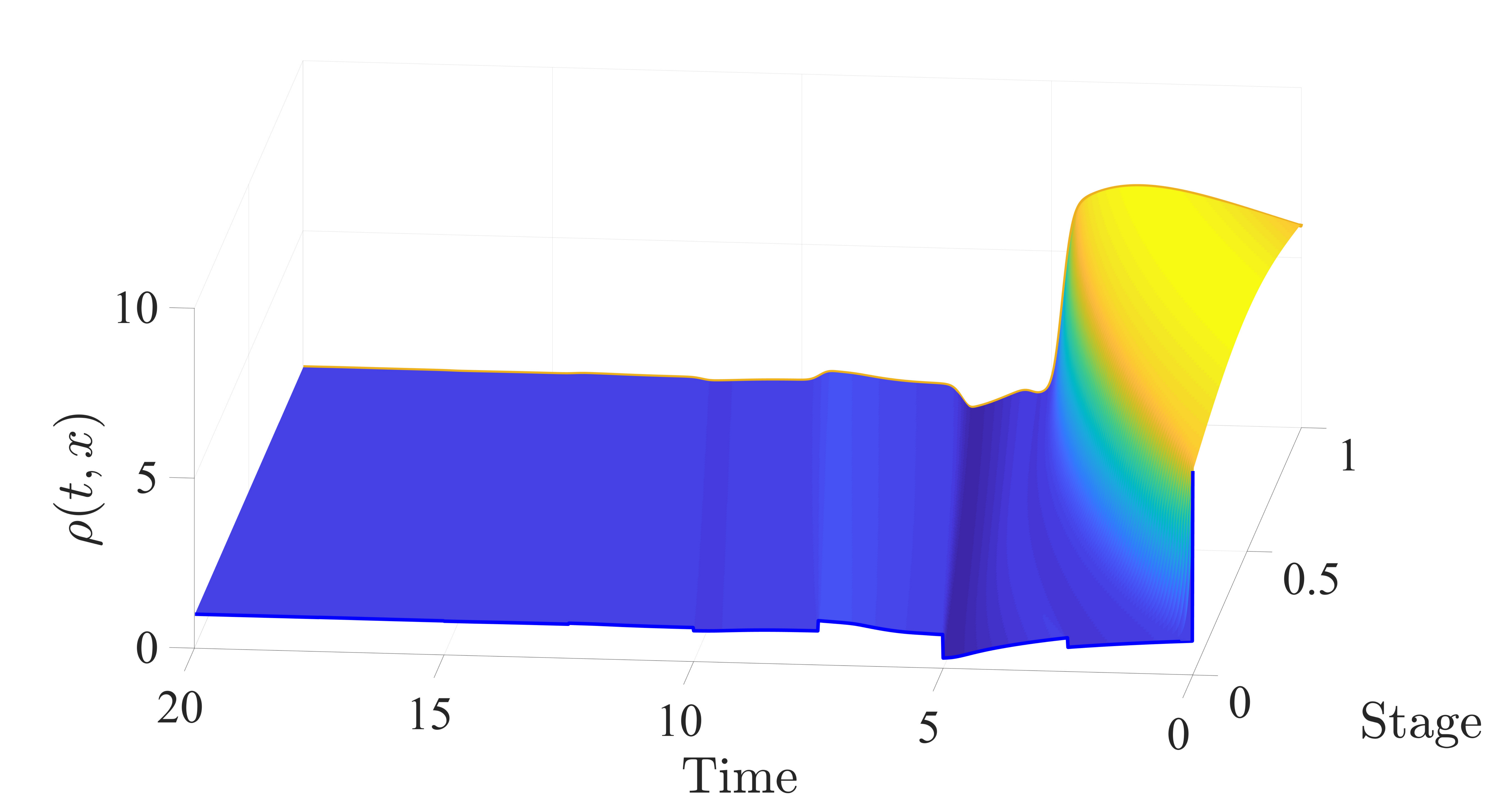}
\caption{The distributed density dynamics with sampling period  $T=1$.}\label{system_figure-sample4} 

\end{figure}

\begin{figure}[htbp]
\centering
\includegraphics[width=0.45\textwidth]{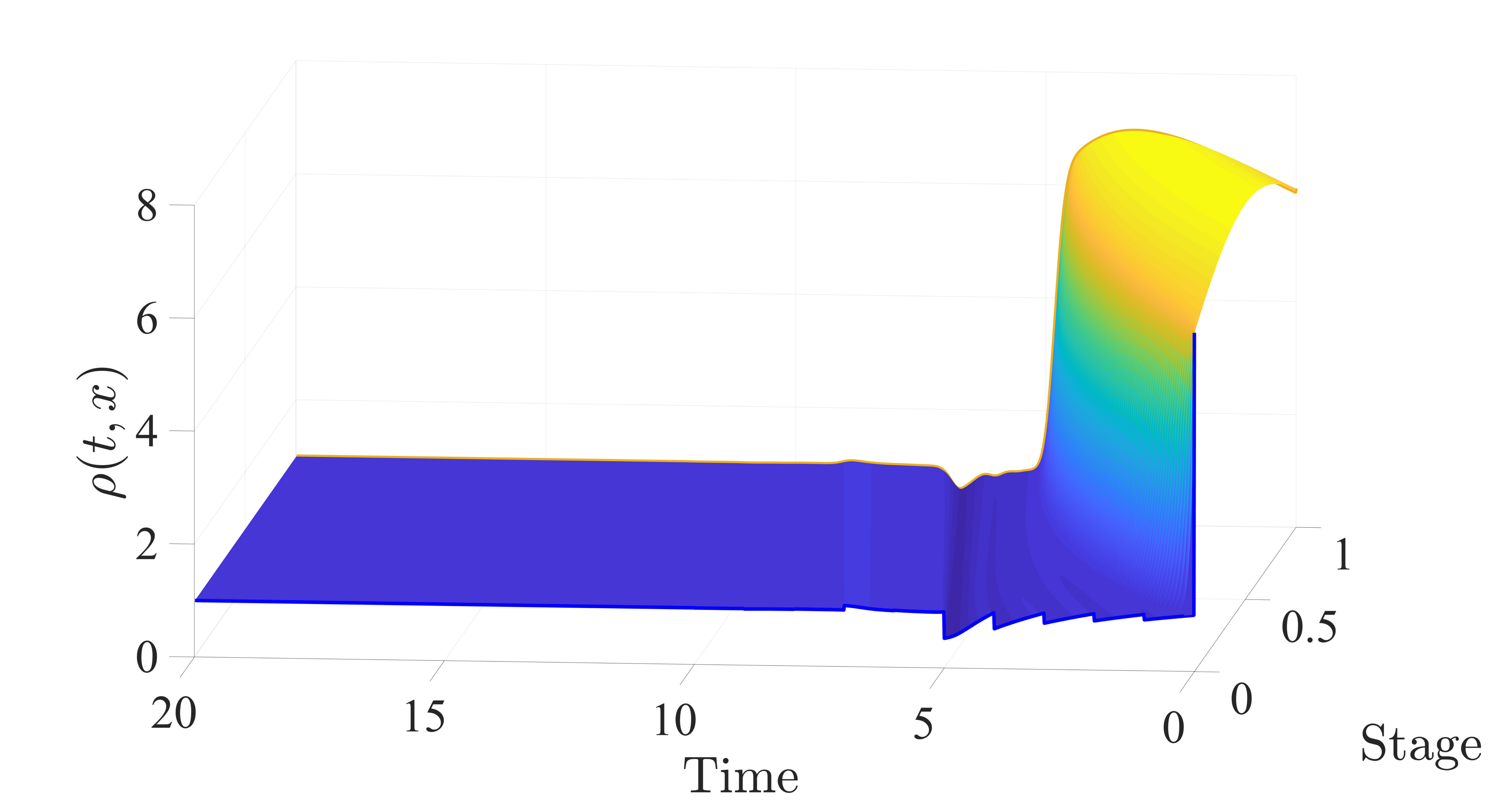}
\caption{The distributed density  dynamics with sampling period  $T=2.5$.}\label{system_figure-sample5} 

\end{figure}

To illustrate the robustness of the control action concerning the sampling schedule, we apply the controller \eqref{GrindEQ__9_} with a periodically updated control action. Here, the simulation is performed under the previous initial condition with an identical function $\lambda$. Two sampling periods are considered, namely,  $T=1$ and $T=2.5$ as shown in Figure \ref{system_figure-sample0}, \textcolor{black}{motivated by the statistics shown in Figure \ref{system_figure6} and Figure \ref{system_figure7}.} The results obtained in Figure \ref{system_figure-sample1} and Figure \ref{system_figure-sample2} prove that the input density converge to the uniform setpoint $\rho_s=1$ and the output flux is also stabilized at the equilibrium. As expected, the $L^2$ norm of the deviation of the state with respect to the uniform equilibrium $\rho_s$ tends to zero (Figure \ref{system_figure-sample3}),
and the distributed density function are stabilized to $\rho_s$ for both the considered sampling periods (Figure \ref{system_figure-sample4} and Figure \ref{system_figure-sample5}). As for the event-triggered control, one can notice that fast sampling ($T=1$) enabled better closed-loop performance.  

\section{Concluding remarks}\label{s6}

In this paper, an even-triggered control algorithm is developed to stabilize the continuum model of a highly re-entrant manufacturing system. The robustness of the proposed controller with respect to the sampling policy is proven to enable the implementation of the classical sampled-data controller with a cyclic update of the control action. Developing an output feedback event-triggered controller for the considered system will be considered in our future works.

\section*{\textcolor{black}{\large Appendix}}

\section{Proof of Theorem \ref{T1}}

 Let $\rho _{0} \in PC^{1} \left([0,1]\right)$, $u>0$ and ${\mathop{\inf }\limits_{x\in (0,1]}} \left(\rho _{0} (x)\right)>0$ be given (arbitrary). Define
\begin{align} \label{GrindEQ__36_} 
\rho _{\max } :={\mathop{\sup }\limits_{x\in (0,1]}} \left(\rho _{0} (x)\right),\quad  \rho _{\min } :={\mathop{\inf }\limits_{x\in (0,1]}} \left(\rho _{0} (x)\right) 
\end{align} 
Consider for each $T\in \left(0,1/\lambda (0)\right]$ the mapping $P:S\to S$, where 
\begin{align} \label{GrindEQ__37_} 
S=\Bigg\{\, W\in C^{0} ([0,T]): 0\le W(t)&\le \rho _{\max } \nonumber\\&+\frac{u}{\lambda (0)}, t\in [0,T]\, \Bigg\} 
\end{align} 
 which maps every $W\in S$ to the function $P(W)\in S$ given by the formula
\begin{align} \label{GrindEQ__38_}
\textcolor{black}{P(W)(t)}=tu+\int _{0}^{a(t)}\rho _{0} (x)dx , \ for \ t\in [0,T]                                             
\end{align} 
 where 
\begin{equation} \label{GrindEQ__39_}
a(t)=1-\int _{0}^{t}\lambda (W(s))ds , \ for \ t\in [0,T] .                                                   
\end{equation}
\textcolor{black}{For all $W,\, V\in S$ for the mapping $P:S\to S$ defined by \eqref{GrindEQ__38_}, \eqref{GrindEQ__39_}, we have that
\begin{equation} \label{GrindEQ__38_0}
|P(W)(t)-P(V)(t)|=\left |      \int_{\bar a(t)}^{ a(t)}\rho_0(x)dx \right|,
\end{equation} 
where
\begin{equation} \label{GrindEQ__39_0}
\bar a(t)=1-\int _{0}^{t}\lambda (V(s))ds , {\ for}  \ t \in [0,T].                                                   
\end{equation}}
From  \eqref{GrindEQ__36_}, \eqref{GrindEQ__39_},  \eqref{GrindEQ__38_0} and \eqref{GrindEQ__39_0}, the following estimate holds
\begin{align} \label{GrindEQ__38_1}
&|P(W)(t)-P(V)(t)|\leq \rho_{\max}\gamma_1(t)\\ 
&\gamma_1(t)=  \left |  \int _{0}^{t} \left[\lambda (W(s))-  \lambda (V(s))\right]ds \right| \nonumber,
\end{align} 
which leads to
\begin{align} \label{GrindEQ__38_2}
&\hspace{-0.5cm}|P(W)(t)-P(V)(t)|\leq t \rho_{\max} \gamma_2\\ 
&\hspace{-0.5cm}\textcolor{black}{\gamma_2=
\sup_{  s\geq 0} \lambda'(s)
\max_{0\leq t \leq T}|W(t)-V(t)|.}
\end{align} 
Hence,  the following inequality holds 
\begin{align} \label{GrindEQ__40_} 
&{\mathop{\max }\limits_{t\in [0,T]}} \left(\left|(P(W))(t)-(P(V))(t)\right|\right)\nonumber\\&\le \rho _{\max } KT{\mathop{\max }\limits_{t\in [0,T]}} \left(\left|W(t)-V(t)\right|\right),
\end{align} 
where $K>0$ is the constant for which $\left|\lambda '(s)\right|\le K$ for all $s\ge 0$. 

Therefore, for $$T=\frac{1}{\lambda (0)+K\rho _{\max } }, $$the mapping $P:S\to S$ is a contraction and Banach's fixed-point theorem implies the existence of a unique $W\in S$ such that 
\begin{align} \label{GrindEQ__42_}
W(t)=tu+\int _{0}^{a(t)}\rho _{0} (x)dx , \ \textrm{for} \ t\in [0,T],                                          
\end{align}
 where $a(t)$ is given by \eqref{GrindEQ__39_}. Notice that definitions \eqref{GrindEQ__36_} and equation \eqref{GrindEQ__42_} as well as the fact that $\lambda (s)\le \lambda (0)$ for all $s\ge 0$, imply the following estimate:
\begin{align} \label{GrindEQ__43_}
W(t)\ge tu+\rho _{\min } (1-\lambda (0)t)>0, \textrm{for} \ t\in [0,T]                                       
\end{align}

\noindent Next define the functions $v:[0,T]\to (0,+\infty )$, $\rho :[0,T]\times [0,1]\to (0,+\infty )$ by means of the equations
\begin{align} \label{GrindEQ__44_}
v(t)=\lambda \left(W(t)\right), \ \textrm{for} \ t\in [0,T] 
\end{align}
\begin{align} \label{GrindEQ__45_} 
\rho (t,x)=\left\{\begin{array}{c} {\rho _{0} \left(x-\int _{0}^{t}v(s)ds \right)\, if\, \int _{0}^{t}v(s)ds<x\leq 1 } \\ {\frac{u}{v\left(\tilde{t}(t,x)\right)} \quad if\quad 0\le x\le \int _{0}^{t}v(s)ds } \end{array}\right.  
\end{align} 
where $\tilde{t}(t,x)\in [0,t]$ is the unique solution of the equation $$x=\int _{\tilde{t}(t,x)}^{t}v(s)ds, $$ for all $$(t,x)\in \Phi :=\left\{\, (t,x)\in [0,T]\times [0,1]\, ,\, x\le \int _{0}^{t}v(s)ds \, \right\}.$$ Notice that $\rho [t]\in PC^{1} ([0,1])$ for each $t\in [0,T]$ and that
\begin{align}\label{GrindEQ__46_}  \int _{0}^{1}\rho (t,x)dx &=\int _{0}^{1-a(t)}\rho (t,x)dx +\int _{1-a(t)}^{1}\rho (t,x)dx\nonumber \\ 
&= u\int _{0}^{1-a(t)}\frac{dx}{v\left(\tilde{t}(t,x)\right)} \nonumber \\
&+\int _{1-a(t)}^{1}\rho _{0} \left(x+a(t)-1\right)dx\nonumber  \\ 
&=ut+\int _{0}^{a(t)}\rho _{0} \left(x\right)dx =W(t) 
\end{align} 
 For the derivation of \eqref{GrindEQ__46_}, we have used formulas \eqref{GrindEQ__39_}, \eqref{GrindEQ__44_}, \eqref{GrindEQ__45_}, \eqref{GrindEQ__42_} and the fact that 
$$\frac{\partial \, \tilde{t}}{\partial \, x} (t,x)=-\frac{1}{v(\tilde{t}(t,x))},$$
 for all
\begin{align*} (t,x)\in \Phi :=\left\{\, (t,x)\in [0,T]\times [0,1]\, ,\, x\le \int _{0}^{t}v(s)ds \, \right\},\end{align*}
with 
\begin{align*}\tilde{t}(t,0)=t \quad  \textrm{and} \quad \tilde{t}\left(t,\int _{0}^{t}v(s)ds \right)=0.\end{align*}

 We next repeat the construction with $\rho _{0} $ replaced by $\rho [T]$. We can construct functions $\tilde{v}:[0,T']\to (0,+\infty )$, $\tilde{\rho }:[0,T']\times [0,1]\to (0,+\infty )$, $\tilde{W}\in C^{0} ([0,T'])$ with $$0\le \tilde{W}(t)\le \left\| \rho [T]\right\| _{\infty } +\frac{u}{\lambda (0)}  \textrm{for} \, t\in [0,T'],$$ and $$T'=\frac{1}{\lambda (0)+K\left\| \rho [T]\right\| _{\infty } } $$ such that
 \begin{align*} \tilde{W}(t)&=\int _{0}^{1}\tilde{\rho }(t,x)dx, \\
  \tilde{v}(t)&=\lambda \left(\tilde{W}(t)\right),\\
   \tilde{W}(t)&=tu+\int _{0}^{\tilde{a}(t)}\rho (T,x)dx, \\
    \tilde{a}(t)&=1-\int _{0}^{t}\lambda (\tilde{W}(s))ds. 
    \end{align*} 
    It is a matter of straightforward calculations to verify that the extensions of $\rho ,v,W$ given by the formulas for $t\in (T,T+T']$, $x\in [0,1]$:
\begin{align*}  
\rho (t,x)&=\tilde{\rho }(t-T,x),\\
 W(t)&=\tilde{W}(t-T), \\
  v(t)&=\tilde{v}(t-T) 
\end{align*} 
 satisfy the following equations for all $t\in [0,T+T']$, $x\in [0,1]$:
\begin{align} \label{GrindEQ__48_} 
v(t)=\lambda \left(W(t)\right) 
\end{align} 
\begin{equation} \label{GrindEQ__49_} 
\rho (t,x)=\left\{\begin{array}{c} {\rho _{0} \left(x-\int _{0}^{t}v(s)ds \right)\quad if\quad \int _{0}^{t}v(s)ds<x\leq 1 } \\ {\frac{u}{v\left(\tilde{t}(t,x)\right)} \quad if\quad 0\le x\le \int _{0}^{t}v(s)ds } \end{array}\right.  
\end{equation} 
\begin{equation} \label{GrindEQ__50_} 
W(t)=\int _{0}^{1}\rho (t,x)dx,  
\end{equation} 
\noindent where $\tilde{t}(t,x)\in [0,t]$ is the unique solution of the equation $$x=\int _{\tilde{t}(t,x)}^{t}v(s)ds, \quad \forall (t,x)\in \Phi,$$ $$  \Phi :=\Bigg\{\, (t,x)\in [0,T+T']\times [0,1]\, ,\\ x\le \int _{0}^{t}v(s)ds \, \Bigg\}.$$

\noindent The construction can be repeated ad infinitum and thus we obtain functions $v:[0,t_{\max } )\to (0,+\infty )$, $\rho :[0,t_{\max } )\times [0,1]\to (0,+\infty )$, $W\in C^{0} ([0,t_{\max } ))$ that satisfy \eqref{GrindEQ__48_}, \eqref{GrindEQ__49_}, \eqref{GrindEQ__50_} for all $t\in [0,t_{\max } )$, $x\in [0,1]$. Moreover, if $t_{\max } <+\infty $ then ${\mathop{\lim \sup }\limits_{t\to t_{\max }^{-} }} \left(\left\| \rho [t]\right\| _{\infty } \right)=+\infty $. Furthermore, $\rho [t]\in PC^{1} ([0,1])$ for each $t\in [0,t_{\max } )$. 

\noindent Finally, formulas \eqref{GrindEQ__48_}, \eqref{GrindEQ__49_} and the facts that $\lambda (s)\le \lambda (0)$ for all $s\ge 0$, ${\mathop{\inf }\limits_{x\in (0,1]}} \left(\rho _{0} (x)\right)>0$ imply that ${\mathop{\inf }\limits_{x\in [0,1]}} \left(\rho (t,x)\right)>0$ for each $t\in [0,t_{\max } )$. Notice that equations \eqref{GrindEQ__2_}, \eqref{GrindEQ__3_}, \eqref{GrindEQ__5_} hold for $t\in [0,t_{\max } )$ with $u(t)\equiv u>0$.

 \textcolor{black}{Expressing  \eqref{GrindEQ__50_} as the integral of \eqref{GrindEQ__49_} over $[0,1]$ and using the change of variables  $\xi=x-\int _{0}^{t}v(s)ds, \, \bar \xi=\tilde t(t,x) $, where $v$ is given by  \eqref{GrindEQ__48_}, we deduce that }  
\begin{align} \label{GrindEQ__51_}
W(t)&=\left(t-\tilde{t}\left(t,\min (1,1-a(t))\right)\right)u\nonumber\\&+\int _{0}^{\max (0,a(t))}\rho _{0} (\xi)d\xi , \ for \ t\in [0,t_{\max } )                          
\end{align}
 where $a(t)$ is given by \eqref{GrindEQ__39_}. 

Equation \eqref{GrindEQ__49_} implies that $\rho$ is not   $C^{1} $ at specific points:
\begin{itemize}
\item The points $(t,x)\in [0,t_{\rm max})\times[0,1]$ for which 
\begin{equation}\label{1a}\xi _{i}=x-\int _{0}^{t }v(s)ds, \end{equation} where $\xi _{i} \in [0,1)$ ($i=0,...,N$) are the points (in increasing order \textcolor{black}{with $\xi_0=0$}) for which $\rho _{0} \in C^{1} \left([0,1]\backslash \{ \xi _{0} ,...,\xi _{N} \} \right)$ due to the lack of regularity of the initial condition $\rho_0(x)$.
\end{itemize}
\begin{itemize}
\item On the other hand, \textcolor{black}{equation} \eqref{GrindEQ__49_}, shows that $\rho$ may not be   $C^{1} $ at the points $(t,x)\in [0,t_{\rm max})\times[0,1]$ for which $\tilde t(t,x)$ is equal to a time where $W$ is not $C^1$. Knowing that \eqref{GrindEQ__51_} can be rewritten as 
\begin{equation} \label{GrindEQ__49_a} 
\hspace{-0.25cm} W (t)=\left\{\begin{array}{c} {ut+ \int _{0}^{a(t)}\rho_0(x) \ dx \quad \textrm{for} \quad a(t)\geq 0} \\ {\left(t-\tilde{t}(t,1)\right)u \quad \textrm{for} \quad a(t)<0} \end{array}\right., 
\end{equation} 
we deduce that the  $W$ is not $C^1$ at the times where $a(t)=0$ and $a(t)$ is equal to a point where $\rho_0$ is discontinuous. Clearly, all discontinuity points of $\rho_0$ are included in the set  $\{ \xi _{0} ,...,\xi _{N} \} $ and \textcolor{black}{therefore the times $\tau \in [0,t_{\rm max})$ }with $a(\tau)\in \{ \xi _{0} ,...,\xi _{N} \}$ are the \textcolor{black}{times} of concern. Consequently, the points 
$(t,x)\in [0,t_{\rm max})\times[0,1]$ for which $\tilde t(t,x)$ is equal to a time where $W$ is not $C^1$ are \textcolor{black}{included in the set of all $(t,x)$ for which}
\begin{align}
a(\tilde t(t,x))\in \{\xi_0,...,\xi_N\}.\end{align}
Since $x=\int_{\tilde t(t,x)}^tv(s)ds$  and $a(t)=1-\int_{0}^1 v(s)ds$, it follows that
\begin{align}
a(\tilde t(t,x))=1 +x -\int_{0}^1 v(s)ds\end{align}
and the discontinuity occurs at the points $(t,x)\in [0,t_{\rm max})\times[0,1]$ satisfying
\begin{align}
\xi_i=1 +x -\int_{0}^1 v(s)ds, \quad \textcolor{black}{i=0,...,N}.\end{align}
Since $$r_i(t)=\xi_i+\int_{0}^1 \lambda (W(s))ds=\xi_i+\int_{0}^1 v(s)ds,$$ 
\textcolor{black}{for $i=0,...,N$, it  follows that $x=r_i(t)-1$.} 
\end{itemize}
Finally, combining the sets defined by \eqref{GrindEQ__6_b} and \eqref{GrindEQ__6_c}, we arrive at \eqref{GrindEQ__6_a}.

\noindent The fact that equation \eqref{GrindEQ__1_} holds for all $(t,x)\in [0,t_{\max } )\times [0,1]\backslash \Omega $ is a direct consequence of formula \eqref{GrindEQ__49_}, the above regularity properties for $\rho ,v,W$ and the fact that $$\frac{\partial \, \tilde{t}}{\partial \, x} (t,x)=-\frac{1}{v(\tilde{t}(t,x))}, \quad \frac{\partial \, \tilde{t}}{\partial \, t} (t,x)=\frac{v(t)}{v(\tilde{t}(t,x))}, $$ for all $(t,x)\in \Phi,$
$$\Phi :=\left\{\, (t,x)\in [0,t_{\max } )\times [0,1]\, ,\, x\le \int _{0}^{t}v(s)ds \, \right\}.$$

\noindent Uniqueness of solution is a consequence of Banach's fixed-point theorem: the fact that equation \eqref{GrindEQ__51_} has a unique solution $W\in C^{0} ([0,t_{\max } ))$. The solution of \eqref{GrindEQ__51_} is constructed step-by-step by using the mapping $P:S\to S$ defined by \eqref{GrindEQ__38_}, \eqref{GrindEQ__39_} and Banach's fixed-point theorem guarantees that $P:S\to S$ has a unique fixed point. 

\noindent The proof is complete.       $\triangleleft $


\begin{thebibliography}{s20}



\bibitem{s1} Yook J.  K., Dawn M. T., and Nandit R.  S. (2002). ``Trading computation for bandwidth: Reducing communication in distributed control systems using state estimators." IEEE transactions on Control Systems Technology, 10 (4), 503-518.

\bibitem{s2} Tabuada P.  (2007). ``Event-triggered real-time scheduling of stabilizing control tasks." IEEE Transactions on Automatic Control, 52 (9), 1680-1685.

\bibitem{s3}  Heemels W. P. M. H.,  Johansson K. H., and Tabuada, P.  (2012). ``An introduction to event-triggered and self-triggered control." IEEE Conference on Decision and Control (CDC), 3270-3285.


\bibitem{s4}  Hespanha J. P., Naghshtabrizi P.,  and  Xu Y. (2007). ``A Survey of Recent Results in Networked Control Systems." Proceedings of the IEEE , 95 (1), 138-162.

\bibitem{s5}  Heemels W. H., Donkers M. C. F., and Teel A. R. (2012). ``Periodic Event-Triggered Control for Linear Systems." IEEE Transactions on Automatic Control, 58 (4), 847-861.


\bibitem{s6}  Peng C., and  Han Q. (2013). ``A Novel Event-Triggered Transmission Scheme and ${L}_{2}$ Control Co-Design for Sampled-Data Control Systems." IEEE Transactions on Automatic Control, 58 (10), 2620-2626.


\bibitem{s7} Seuret A., and Christophe P. (2011). ``Event-triggered sampling algorithms based on a Lyapunov function."  IEEE Conference on Decision and Control and European Control Conference,   6128-6133.

\bibitem{s8}  Aström K.  J., and Bernhardsson B. P. (1999). ``Comparison of periodic and event based sampling for first-order stochastic systems."  IFAC World Congress, 32 (2)   5006-5011.

\bibitem{s9}  Arz\'en K. E. (1999). ``A simple event-based pid controller."  IFAC World Congress, 32 (2),   423-428.

 
\bibitem{s10} Tallapragada P.,  and  Chopra N. (2013). ``On Event Triggered Tracking for Nonlinear Systems."   IEEE Transactions on Automatic Control, 58 (9), 2343-2348.

\bibitem{s11}Postoyan R.,  Tabuada P., Nesi\'c  D., and  Anta A. (2015). ``A Framework for the Event-Triggered Stabilization of Nonlinear Systems."   IEEE Transactions on Automatic Control, 60 (4), 982-996.

\bibitem{s12}Abdelrahim M.,  Postoyan R.,  Daafouz J., and  Nesi\'c D. (2016). ``Stabilization of Nonlinear Systems Using Event-Triggered Output Feedback Controllers."  IEEE Transactions on Automatic Control,  61 (9), 2682-2687.

\bibitem{s12}Abdelrahim M.,  Postoyan R.,  Daafouz J., and  Nesi\'c D. (2016). ``Stabilization of Nonlinear Systems Using Event-Triggered Output Feedback Controllers."  IEEE Transactions on Automatic Control,  61 (9), 2682-2687.


\bibitem{s13} Hetel L.,  Fiter C., Omran H., Seuret A.,  Fridman E.,  Richard J-.P., and   Niculescu S. L. (2017). ``Recent developments on the stability of systems with aperiodic sampling: An overview." Automatica, 76, 309-335.


\bibitem{s14}Borgers D. P., and  Heemels W. P. M. H. (2014). ``Event-separation properties of event-triggered control systems."  IEEE Transactions on Automatic Control, 59 (10), 2644-2656.


 \bibitem{s15} Nowzari C.,  Garcia E., and Cort\'es J. (2019).
``Event-triggered communication and control of networked systems for multi-agent consensus." Automatica, 105, 1-27.

\bibitem{s16} Xing L.,  Wen C.,  Liu Z.,  Su H., and Cai J. (2017) ``Event-Triggered Adaptive Control for a Class of Uncertain Nonlinear Systems." IEEE Transactions on Automatic Control, 62 (4), 2071-2076.

\bibitem{s34}Blevins T.,  Nixon M., and  Wojsznis W. (2015). ``Event Based Control Applied to Wireless Throttling Valves."  International Conference on Event-Based Control, Commmunication, and Signal Processing, 1-6



\bibitem{s36}Pawlowski A.,  Guzm\`an J. L.,  Berenguel M., and Dormido S. (2016). ``Event-based Generalized Predictive Control."  Event-based Control and Signal Processing,   Boca Raton, FL: CRC Press, 151?176.

 \bibitem{s37}Guerrero-Castellanos J. F., Vega-Alonzo A.,  Marchand N., Durand S. ,  Linares-Flores J., and  Mino-Aguilar G. (2017).``Real-time event-based formation control of a group of VTOL-UAVs." International Conference on Event-Based Control, Communication and Signal Processing, 1-8.
 
 \bibitem{s35}Boisseau B.,  Durand S., Martinez-Molina J. J.,  Raharijaona T., and N. Marchand. (2015). ``Attitude Control of a Gyroscope Actuator Using Event-based Discrete-time Approach."  International Conference on Event-Based Control, Communication, and Signal Processing, 1-6.


\bibitem{s17} Anton S., and  Fridman E. (2015). ``Distributed event-triggered control of transport-reaction systems." IFAC Conference on Modelling, Identification and Control of Nonlinear Systems (MICNON), 48 (11), 593-597.

\bibitem{s18}Anton S., and  Fridman E. (2016). ``Distributed event-triggered control of diffusion semilinear PDEs." Automatica,  68,  344-351.

 \bibitem{s29} Wen K, and  Fridman E. (2018). ``Distributed sampled-data control of Kuramoto Sivashinsky equation." Automatica, 95, 514-524.


\bibitem{s19} Davo M. A., Bresch-Pietri D., Prieur C., and Di Meglio F. (2018). ``Stability Analysis of a ${\text {2}}\times {\text {2}} $ Linear Hyperbolic System With a Sampled-Data Controller via Backstepping Method and Looped-Functionals." IEEE Transactions on Automatic Control 64 (4), 1718-1725.

 
  
 
\bibitem{s20} Wang, J.-W. (2019) ``Observer-based boundary control of semi-linear parabolic PDEs with non-collocated distributed event-triggered observation." Journal of the Franklin Institute 356 (17),  10405-10420.
 
 \bibitem{s32} Wang J.-W., and  Wang J.-M. (2019) ``Mixed $H^2/H^\infty$ sampled-data output feedback control design for a semi-linear parabolic PDE in the sense of spatial $L^\infty$ norm." Automatica, 103,  282-293.
 
 
\bibitem{s22} Espitia N., Girard A., Marchand N., and Prieur C. (2017). ``Event-based boundary control of a linear $2\times 2$ hyperbolic system via backstepping approach." IEEE Transactions on Automatic Control 63 (8), 2686-2693.
 
 
 
 
 \bibitem{s23}Espitia N. (2020) ``Observer-based event-triggered boundary control of a linear $2\times 2$ hyperbolic systems." Systems \& Control Letters 138,104668.
 
 
 
\bibitem{s24} Espitia N., Aneel T., and \textcolor{black}{Tarbouriech S.}  (2017). ``Stabilization of boundary controlled hyperbolic PDEs via Lyapunov-based event triggered sampling and quantization."  IEEE  Conference on Decision and Control (CDC), 1266-1271.
 
 \bibitem{s25}Espitia N., Girard A., Marchand N., and Prieur C. (2016) ``Event-based control of linear hyperbolic systems of conservation laws." Automatica, 70, 275-287.
 
 
 \bibitem{s25a} \textcolor{black}{Espitia, N., Karafyllis, I., and  Krstic, M. (2019). ``Event-triggered
boundary control of constant-parameter reaction-diffusion
PDEs: a small-gain approach." arXiv:1909.10472.}
 
 \bibitem{s26} Baudouin L., Marx S., and Tarbouriech S. (2019). ``Event-triggered damping of a linear wave equation."  IFAC Workshop on Control of Systems Governed by Partial Differential Equations CPDE, 52 (2), 58-63.
 
 
 
 
 
 
 \bibitem{s27} Karafyllis I., and Krstic M. (2018). ``Sampled-data boundary feedback control of 1-D parabolic PDEs." Automatica, 87, 226-237.
 
 
  \bibitem{s28}Karafyllis I., and Krstic M. (2017).  ``Sampled-data boundary feedback control of 1-D linear transport PDEs with non-local terms." Systems \& Control Letters, 107, 68-75.
 
 
 \bibitem{s33} Karafyllis I., and Krstic M. (2020). ``Stability results for the continuity equation." Systems \& Control Letters, 135,104594.
 
  \bibitem{s33x}Karafyllis I. and  Krstic M. (2017). ``Stability of Integral Delay Equations and Stabilization of Age Structured Models." ESAIM Control, Optimisation and Calculus of Variations, 23,  1667-1714
 
 
 
 
  \bibitem{s30}Ahmed-Ali T., Karafyllis I., Giri F., Krstic M.,  and  Lamnabhi-Lagarrigue, F. (2017).  "Exponential stability analysis of sampled-data ODE-PDE systems and application to observer design." IEEE Transactions on Automatic Control 62 (6), 3091-3098.
 
 \bibitem{s31}  Zhiyuan Y., and  El-Farray N. H. (2012). ``A predictor-corrector approach for multi-rate sampled-data control of spatially distributed systems."  IEEE Conference on Decision and Control (CDC), 2908-2913.
 
 \bibitem{s21}  Da X., and  El-Farray N. H. (2017). ``Resource-aware fault accommodation in spatially-distributed processes with sampled-data networked control systems."  IEEE American Control Conference (ACC), 1809-1814.

 
 \bibitem{s38} La Marca M., Armbruster  D., Herty  M., and  Ringhofer  C. (2010). ``Control of continuum models of production systems." IEEE Transactions on Automatic Control, 55 (11), 2511-2526.

\bibitem{s39} D'Apice C., Kogut P. I., and Manzo R. (2016). On optimization of a highly re-entrant production system. Networks \& Heterogeneous Media, 11(3), 415-445.

\bibitem{s40} Shang P., and Wang Z. (2011).`` Analysis and control of a scalar conservation law modeling a highly re-entrant manufacturing system." Journal of Differential Equations, 250 (2), 949-982.

\bibitem{s41} Coron J. M., and Wang  Z. (2012). ``Controllability for a scalar conservation law with nonlocal velocity." Journal of Differential Equations, 252 (1), 181-201.

\bibitem{s42} Coron, J. M., and  Wang, Z. (2013). ``Output feedback stabilization for a scalar conservation law with a nonlocal velocity." SIAM Journal on Mathematical Analysis, 45 (5), 2646-2665.


\bibitem{s43} Armbruster D., Gottlich S., and  Herty M. (2011). ``A scalar conservation law with discontinuous flux for supply chains with finite buffers." SIAM Journal on Applied Mathematics, 71 (4), 1070-1087.

\bibitem{s44} Armbruster D., Marthaler D., and Ringhofer C. (2003). ``Kinetic and fluid model hierarchies for supply chains." Multiscale Modeling \& Simulation, 2(1), 43-61.


\bibitem{s45} Xu X., Ni D., Yuan Y., and Dubljevic, S. (2018). ``PI-control design of continuum models of production systems governed by scalar hyperbolic partial differential equation." IFAC Symposium on Advanced Control of Chemical Processes ADCHEM  51(18), 584-589.

\bibitem{s46} Diagne M., Bekiaris-Liberis N., and Krstic M. (2017). ``Compensation of input delay that depends on delayed input."  Automatica, 85, 362-373.
 
  
  \bibitem{s47}
  Zhou M. and  DiCesare F. (2012). ``Petri net synthesis for discrete event control of manufacturing systems."  Springer Science \& Business Media, 204, 1-201.

\bibitem{s48} Hu H., and  Zhou M. (2014). ``A Petri net-based discrete-event control of automated manufacturing systems with assembly operations." IEEE Transactions on Control Systems Technology, 23, 2, 513-524.


\bibitem{s49} Lefeber E., Van Den Berg R. A., and Rooda J. E. (2004). ``Modeling, validation and control of manufacturing systems." IEEE  American Control Conference, 5,  4583-4588. 



\bibitem{s50} Graves S. C. (1986). ``A tactical planning model for a job shop." Operations Research, 34, 4, 522-533.


\bibitem{s51}Jacobs J. H., Etman L. F. P., Van Campen E. J. J., and Rooda, J. E. (2003). ``Characterization of operational time variability using effective process times." IEEE Transactions on Semiconductor Manufacturing, 16, 3, 511-520.

\end{thebibliography}
\end{document}